\documentclass[amscd,amssymb,11pt,reqno]{amsart}

\newcommand {\eps}{\varepsilon}

\newtheorem{Thm}{Theorem}[section]

\newtheorem{Lem}[Thm]{Lemma}
\newtheorem{Prop}[Thm]{Proposition}

\newtheorem{Rmk}[Thm]{Remark}

\usepackage{graphicx}
\usepackage{latexsym}
\begin{document}

%\begin{Large}

\vspace{1.5 cm}

\title[Busemann's intersection inequality]
      {Busemann's intersection inequality  in hyperbolic and spherical spaces}

\author{Susanna Dann, Jaegil Kim, and Vladyslav Yaskin}\thanks{The second and third named authors are supported in part by NSERC}

\address{Vienna University of Technology, Wiedner Hauptstrasse 8-10, 1040 Vienna, Austria} \email{susanna.dann@tuwien.ac.at}

\address{Department of Mathematical and Statistical Sciences, University of Alberta, Edmonton, Alberta T6G 2G1, Canada} \email{jaegil@ualberta.ca}

\address{Department of Mathematical and Statistical Sciences, University of Alberta, Edmonton, Alberta T6G 2G1, Canada} \email{vladyaskin@math.ualberta.ca}

\date{\today}


\begin{abstract}
Busemann's intersection inequality asserts that the only maximizers of the integral
$\int_{S^{n-1}} |K\cap \xi^\perp|^n d\xi$ among all convex bodies of a fixed volume in $\mathbb R^n$ are centered ellipsoids. We study this question in the hyperbolic and spherical spaces, as well as general measure spaces.

\end{abstract}

\maketitle

\section{Introduction}

Let $K$ be a convex body in $\mathbb R^n$ that contains the origin in its interior. The following is known as Busemann's intersection inequality:
\begin{equation}\label{BI}
\int_{S^{n-1}} |K\cap \xi^\perp|^n \ d\xi \le c_n |K|^{n-1},
\end{equation}
with equality if and only if $K$ is a centered ellipsoid; see \cite{Bu}. Here, $c_n=n \kappa_{n-1}^n/\kappa_n^{n-2}$,  where
$\kappa_n$ is the volume of the unit Euclidean  ball
$B^n$   in $\mathbb R^n$, and $|A|$ stands for the volume (in the appropriate dimension) of a set $A$.

In fact, the inequality (\ref{BI}) is true for a larger class of sets, in particular, star bodies; see \cite[p.373]{Ga}. In a slightly different form, (\ref{BI}) can be stated as follows.
Centered ellipsoids in $\mathbb R^n$ are the only maximizers of the quantity
\begin{equation}\label{int}
\int_{S^{n-1}} |K\cap \xi^\perp|^n d\xi
\end{equation}
 in the class of star bodies of a fixed volume.
In this paper we study this question in the hyperbolic space $\mathbb
H^n$ and the sphere $\mathbb S^n$ (or, more precisely, a hemisphere
$\mathbb S^n_+$, as explained in the next section). We show that in
$\mathbb H^n$  centered  balls are the unique maximizers of
(\ref{int}) in the class of star bodies of a fixed volume. On the
sphere the situation is different. In $\mathbb S^2_+$   centered
balls are in
fact the unique minimizers (in the class of origin-symmetric star bodies). The maximizers of
(\ref{int}) in the class of origin-symmetric star bodies in $\mathbb
S^2_+$ are cones (see Section 2 for the definition).  The
maximizers of (\ref{int}) in the class of origin-symmetric convex
bodies in $\mathbb S^2_+$ are   lunes. It is surprising that in
$\mathbb S^n_+$ with $n\ge 3$   centered balls are neither maximizers nor minimizers,
even in the class of origin-symmetric convex bodies. We also obtain
 an optimal lower bound for (\ref{int}) in the class of star bodies in $\mathbb S^n_+$, $n\ge 3$, of a given volume and describe the
 equality cases.
Finally, we prove a version of Busemann's intersection inequality (together with the
equality cases)
for general measures on $\mathbb R^n$ and $\mathbb H^n$. An important
special  case is that of the Gaussian measure on $\mathbb R^n$.

For the history of Busemann's inequality, its applications, and recent developments the reader is referred to  \cite{DPP}, \cite{Ga},  \cite{Ga1}, \cite{Gri}, \cite{Lu}.

It is interesting to note that, in the context of the Busemann-Petty problem, the sphere and the Euclidean space are similar in the sense that the positive answer holds in the same dimensions, while the hyperbolic space exhibits a different behaviour; see \cite{Y}. For Busemann's intersection inequality, the hyperbolic space is similar to the Euclidean space, while the sphere is not.

\section{Preliminaries}

We will start with some basic facts about the hyperbolic space $\mathbb H^n$ and the sphere $\mathbb S^n$. In $\mathbb H^n$ any two points are connected by a unique geodesic line. Thus we say that a set $A$ in $\mathbb H^n$ is {\em convex} if for any two points that belong to $A$, the geodesic segment that connects these two points also belongs to $A$.
On the sphere $\mathbb S^n$ the uniqueness property of geodesics does not hold. To overcome this obstacle, one can consider an open hemisphere, where geodesic convexity is well-defined. We will however need a larger class of convex sets. Let $\mathbb S^n_+$ be the closed upper hemisphere of  the unit sphere in $\mathbb R^{n+1}$. A set $A\subset \mathbb S^n_+$ is said to be {\em convex} if $A$ is obtained as the intersection of the hemisphere $\mathbb S^n_+\subset \mathbb R^{n+1}$ and  a convex cone in $\mathbb R^{n+1}$ with vertex at the origin. One can see that the class of convex sets contains all geodesically convex sets in the interior of  $\mathbb S^n_+$, as well as certain sets that have common points with the relative boundary of $\mathbb S^n_+$.

To treat both spaces simultaneously, we will denote them by $\mathbb M_\delta^n$, where $\delta = -1$ corresponds to $\mathbb H^n$  and $\delta = 1$ to $\mathbb S^n_+$.
The definition of star bodies in $\mathbb M_\delta^n$ is similar to that in the Euclidean space. More precisely, first fix an origin $o$ in $\mathbb M_\delta^n$ and denote by $T_o \mathbb M_\delta^n$ the tangent space to $\mathbb M_\delta^n$ at $o$. In the case of $\mathbb S^n_+$, the origin is the center of the   hemisphere. We say that a subset $K$ of $\mathbb M_\delta^n$ is {\em star-shaped} (with respect to $o$), if every geodesic passing through $o$ intersects $K$ in a connected segment containing $o$.  Consider the unit sphere $S^{n-1}$ in $T_o \mathbb M_\delta^n$, i.e. the set of vectors in $T_o\mathbb M_\delta^n$ that have length one with respect to the metric of $\mathbb M_\delta^n$. For each $\xi\in S^{n-1}\subset T_o \mathbb M_\delta^n$, consider the geodesic ray $\gamma$ emanating from $o$ with tangent vector $\xi$. The {\it radial function} of $K$ in the direction of $\xi$ is defined by
$$\rho_K(\xi) = \sup_{x\in \gamma\cap K} \mathrm{d}(x,o),$$
where $\mathrm{d}$ is the metric in $\mathbb M_\delta^n$.

We say that $K$ is a {\it star body} in $\mathbb M_\delta^n$ if $K$ is compact, star-shaped with respect to $o$, and its radial function is positive and continuous. A star-shaped set $K$ is {\it origin-symmetric} if  $\rho_K(\xi) = \rho_K(-\xi)$ for all $\xi\in S^{n-1}$.
For each $\xi\in S^{n-1}\subset T_o\mathbb M_\delta^n$,  we  denote by $\xi^\perp$ the  (unique) totally geodesic submanifold of $\mathbb M_\delta^n$ passing through $o$, whose normal vector at $o$ is $\xi$. As in the Euclidean space, we will often refer to $\xi^\perp$ as a hyperplane.

    The space $\mathbb M_\delta^n$ for both $\delta = \pm 1$ can be identified with   the  unit ball $B^n\subset \mathbb{R}^n$  (open when $\delta=-1$ and closed when $\delta = 1$) equipped with the metric
\begin{eqnarray}\label{eqn:metric}
ds^2=4 \frac{dx_1^2+\cdots +dx_n^2}{(1+\delta\ (x_1^2+\cdots+x_n^2))^{2}}.
\end{eqnarray}

In this model star bodies in $\mathbb M_\delta^n$ correspond to Euclidean star bodies in the ball $B^n$.
Thinking of $K$ as a body in $B^n\subset \mathbb R^n$, we can use standard Euclidean concepts.
The {\it Minkowski functional} of $K\subset \mathbb R^n$ is defined by $$\|x\|_K=\min\{ a\ge 0: x\in aK\}, \quad x\in \mathbb R^n.$$
In  metric (\ref{eqn:metric}) the volume of $K$ equals
\begin{eqnarray*}
\mathrm{vol} (K)
 =2^{n} \int_{S^{n-1} }
            \int_0^{\|\theta\|_K^{-1}} \frac{r^{n-1} dr}{(1+\delta\, r^2)^{n}}\, d\theta,
\end{eqnarray*}
where, once again, $\|\theta\|_K$ is well-defined, since $K$ is a body in $B^d\subset \mathbb R^n$.

The volume of the section of $K$ by the hyperplane $\xi^\perp$ is given by
\begin{eqnarray*}
\mathrm{vol} (K\cap \xi^\perp)
 =2^{n-1} \int_{S^{n-1}\cap\xi^\perp}
            \int_0^{\|\theta\|_K^{-1}} \frac{r^{n-2} dr}{(1+\delta\, r^2)^{n-1}}\, d\theta.
\end{eqnarray*}

When   a body $K$  lies in the Euclidean space, by $|K|$ and $|K\cap \xi^\perp|$ we denote its volume and the volume of its sections with respect to the Euclidean metric.

We will also use another coordinate system for $\mathbb M_\delta^n$, in which
\begin{eqnarray*}
\mathrm{vol} (K)
 =  \int_{S^{n-1} }
            \int_0^{\rho_K(\theta)} s_\delta^{n-1} (r)\, dr\, d\theta,
\end{eqnarray*}
and
\begin{eqnarray*}
\mathrm{vol} (K\cap \xi^\perp)
 =  \int_{S^{n-1}\cap\xi^\perp}
            \int_0^{\rho_K(\theta)} s_\delta^{n-2} (r)\, dr\,  d\theta,
\end{eqnarray*}
where $s_1(r)=\sin r$ and $s_{-1}(r)=\sinh r$.

Below we introduce some concepts specific to the spherical case.
We say that a star-shaped subset $K$ of $\mathbb S^n_+$ is a {\em spherical cone} (or simply a {\em cone} in $\mathbb{S}^n_+$) if the radial function of $K$ is   equal to $\pi/2$ on its support.
In this case, the   support of the radial function, which is a subset of $S^{n-1}$,  is called the   {\em base} of the cone.

Two geodesic lines in $\mathbb{S}^2_+$ are said to be {\em parallel}
to each other if they never meet in the interior of $\mathbb{S}^2_+$,
i.e., if they meet only on the boundary  of $\mathbb{S}^2_+$. A region in
$\mathbb{S}^2_+$ is called a {\em lune}   if its boundary consists of two parallel geodesic lines.

We will use some tools of harmonic analysis, in particular, the spherical Radon transform and
spherical harmonics. Recall that for a function $f\in C(S^{n-1})$ its
spherical Radon transform $Rf$ is the function on $S^{n-1}$ given by
$$Rf(\xi) = \int_{S^{n-1}\cap \xi^\perp} f(\theta)\, d\theta, \qquad
\xi \in S^{n-1}.$$

On the spherical harmonics of a fixed degree the operator $R: C(S^{n-1})
\to C(S^{n-1})$ acts as a multiple of the identity. Namely, if  $H_k$
is a spherical harmonic of an even degree $k$, then
$$RH_k = \lambda_k H_k,$$
where
$$\lambda_k  =
\frac{(-1)^{k/2} 2\pi^{(n-2)/2} \Gamma((k+1)/2)}{ \Gamma((n+k-1)/2)};$$
cf., for example, \cite[p.103]{Gr}. For odd $k$ all the multipliers are equal to zero.

Observe that $|\lambda_k|$, for even $k$, form a strictly decreasing sequence tending
to zero as  $k\to\infty$. Since $\lambda_0=|S^{n-2}|$, the latter implies
that
\begin{equation}\label{L2bound}
\|Rf\|_{L^2(S^{n\!-\!1})}\le |S^{n-2}|\  \|f\|_{L^2(S^{n\!-\!1})}.
\end{equation}
We will   often use the following relation:
\begin{equation}\label{self-adj}
\int_{S^{n-1}}  \int_{S^{n-1}\cap\xi^\perp} g(\theta )\, d\theta \,
d\xi = |S^{n-2}| \int_{S^{n-1}}    g(\xi )\, d\xi.
\end{equation}

\section{Busemann's intersection inequality in the hyperbolic space}

Before we state our results, let us introduce the following functions.
For $n\ge 1$, let the function $F_n:[0,1)\to [0,\infty)$ be defined by $$F_n(t) = \int_0^t \frac{r^{n-1}}{(1-r^2)^n} dr.$$
If $n\ge 2$, consider the function
$$G(t) = \left[F_{n-1}\left(F_n^{-1}(t)\right) \right]^{n/(n-1)},\quad t\in[0,\infty).$$

\begin{Lem}\label{G_concave}
The function $G$ is concave.
\end{Lem}
\proof Since $$G(F_n (t)) = \left[F_{n-1}\left(t\right) \right]^{n/(n-1)},\quad t\in[0,1),$$ we have
$$G^{'}(F_n(t)) \frac{t^{n-1}}{(1-t^2)^n} = \frac{n}{n-1} \left[F_{n-1}\left(t\right) \right]^{1/(n-1)}\frac{t^{n-2}}{(1-t^2)^{n-1}}.$$
Thus, we obtain
$$G^{'}(F_n(t))   = \frac{n}{n-1} \left[F_{n-1}\left(t\right) \right]^{1/(n-1)}\left(\frac{1}{t}-t\right).$$
Differentiating one more time, we get
\begin{multline} G^{''}(F_n(t))  \frac{t^{n-1}}{(1-t^2)^n}  = \frac{n}{(n-1)^2} \left[F_{n-1}\left(t\right) \right]^{(-n+2)/(n-1)}\frac{t^{n-3}}{(1-t^2)^{n-2}}\\ +\frac{n}{n-1} \left[F_{n-1}\left(t\right)
\right]^{1/(n-1)}\left(-\frac{1}{t^2}-1\right).
\end{multline}
To show that $G^{''}<0$, we need to prove the inequality $$(n-1)F_{n-1}(t) \ge \frac{t^{n-1}}{(1-t^2)^{n-2}(1+t^2)}, \quad t\in [0,1).$$
The functions on both sides of the inequality vanish at $t=0$.  Thus it is enough to show that the same inequality holds for their derivatives,
i.e., we need to show that
\begin{multline}
(n-1)\frac{t^{n-2}}{(1-t^2)^{n-1}} \ge  (n-1)t^{n-2} (1-t^2)^{-n+2}(1+t^2)^{-1} \\ - 2 t^{n} (-n+2) (1-t^2)^{-n+1}(1+t^2)^{-1} - 2 t^{n} (1-t^2)^{-n+2}(1+t^2)^{-2},
\end{multline}
or equivalently,
\begin{eqnarray*}
(n-1)(1+t^2)^{2}  & \ge&  (n-1)  (1-t^2) (1+t^2)  \\
& &  - 2 t^{2} (-n+2)   (1+t^2)  - 2 t^{2} (1-t^2)  \\
%$$ = (n-1) - (n-1)t^4 +  (2n-4) t^2 +  (2n-4)t^4 - 2t^2 + 2 t^4$$
& = &(n-1)  +  (2n-6) t^2 +  ( n-1)t^4  .
\end{eqnarray*}
This is true, since $$(n-1)  +  (2n-6) t^2 +  ( n-1)t^4 \le (n-1)  +  2( n-1) t^2 +  ( n-1)t^4  = (n-1) (1+t^2)^2.$$
\qed

Let us introduce a new function
$$H(t) = \left[  G\left( \frac{t}{2^n |S^{n-1}|}   \right)\right]^{n -1}, \quad t\ge 0.$$

\begin{Thm}
Let $K$ be a star body in $\mathbb H^n$ for $n\ge2$. Then
$$\int_{S^{n-1}}  \mathrm{vol} (K\cap \xi^\perp) ^n d\xi\le C_n H(\mathrm{vol}(K)) ,$$
with equality if and only if $K$ is a ball centered at the origin. Here, $$C_n=  |S^{n-1}|^{n-1} n^2 2^{n(n-1)} \left(1-\frac1n\right)^n \frac{\kappa_{n-1}^n}{\kappa_n^{n-2}}.$$

\end{Thm}
\proof
We will work in the Poincar\'e model of $\mathbb H^n$ in the Euclidean ball $B^n$, so
$$\int_{S^{n-1}} \mathrm{vol}(K\cap \xi^\perp)^n d\xi= \int_{S^{n-1}} \left(\int_{S^{n-1} \cap \xi^\perp} \int_0^{\|\theta\|_K^{-1}} \frac{2^{n-1}  r^{n-2} }{(1-r^2)^{n-1}} dr\, d\theta \right)^n d\xi.$$
Let $L$ be a  star body in $\mathbb R^n$ whose Minkowski functional is given by
$$ \|\theta\|_L^{-1} =  \left( \int_0^{\|\theta\|_K^{-1}} \frac{2^{n-1}  r^{n-2} }{(1-r^2)^{n-1}} dr\right)^{1/(n-1)}.$$
Then
\begin{align*} \int_{S^{n-1}}  \mathrm{vol}(K\cap \xi^\perp) ^n d\xi & = \int_{S^{n-1}} \left(\int_{S^{n-1} \cap \xi^\perp} \|\theta\|_L^{-n+1}\, d\theta \right)^n d\xi\\
 & = (n-1)^n \int_{S^{n-1}}  |L \cap \xi^\perp|  ^n d\xi.
 \end{align*}
We can now use Busemann's intersection inequality in $\mathbb R^n$.
\begin{align*}& \le n (n-1)^n \frac{\kappa_{n-1}^n}{\kappa_n^{n-2}}  |L|^{n -1}= n (n-1)^n \frac{\kappa_{n-1}^n}{\kappa_n^{n-2}} \left( \frac1n \int_{S^{n-1}}  \|\theta\|_L^{-n} d\theta \right)^{n -1}\\
& = n^2 2^{n(n-1)} \Big(1-\frac1n\Big)^n \frac{\kappa_{n-1}^n}{\kappa_n^{n-2}} \left[ \int_{S^{n-1}} \left( \int_0^{\|\theta\|_K^{-1}} \frac{ r^{n-2} \,dr }{(1-r^2)^{n-1}}\right)^{\frac{n}{n-1}} d\theta \right]^{n -1},
\end{align*}
where the equality case in the above inequality  holds if and only if $L$ is a centered ellipsoid.
\begin{align*}&= n^2 2^{n(n-1)} \Big(1-\frac1n\Big)^n \frac{\kappa_{n-1}^n}{\kappa_n^{n-2}} \left[  \int_{S^{n-1}} \left( F_{n-1}(\|\theta\|_K^{-1}) \right)^{n/(n-1)} d\theta \right]^{n -1}\\
 &= n^2 2^{n(n-1)} \Big(1-\frac1n\Big)^n \frac{\kappa_{n-1}^n}{\kappa_n^{n-2}} \left[  \int_{S^{n-1}} G \left( F_{n}(\|\theta\|_K^{-1}) \right)  d\theta \right]^{n -1}\\
&=|S^{n-1}|^{n-1} n^2 2^{n(n-1)} \Big(1-\frac1n\Big)^n \frac{\kappa_{n-1}^n}{\kappa_n^{n-2}} \left[   \int_{S^{n-1}} G \left( F_{n}(\|\theta\|_K^{-1}) \right)  \frac{d\theta}{|S^{n-1}|} \right]^{n -1}.
\end{align*}
By Jensen's inequality,
$$\le |S^{n-1}|^{n-1} n^2 2^{n(n-1)} \Big(1-\frac1n\Big)^n \frac{\kappa_{n-1}^n}{\kappa_n^{n-2}} \left[  G\left( \int_{S^{n-1}}   F_{n}(\|\theta\|_K^{-1})  \frac{d\theta}{|S^{n-1}|} \right)\right]^{n -1}.$$

Since $G$ is not affine (the proof of Lemma \ref{G_concave} shows that its second derivative is strictly negative), the equality holds if and only if $F_{n}(\|\theta\|_K^{-1})$ is constant on the sphere, i.e.
$\|\theta\|_K^{-1}$ is constant on the sphere, meaning that $K$ is a ball centered at the origin.

Finally,
\begin{align*}\int_{S^{n-1}} \mathrm{vol}(K\cap \xi^\perp)^n d\xi& \le C_n \left[  G\left( \int_{S^{n-1}}   \int_0^{\|\theta\|_K^{-1}} \frac{r^{n-1}\,dr}{(1-r^2)^n} \frac{d\theta}{|S^{n-1}|}
\right)\right]^{n -1}\\
&=  C_n \left[  G\left( \frac{\mathrm{vol}(K)}{2^n |S^{n-1}|} \right)\right]^{n -1} = C_n H(\mathrm{vol}(K) ).
\end{align*}

\qed

 \section{Busemann's intersection inequality in the spherical space: first glance}

The approach used in the previous section does not work for the spherical space, because the analogue of the  function  $G$ in the spherical setting fails to possess the needed convexity properties. However, it is possible to obtain weaker estimates.

Let us introduce a function $F$ by the formula
$$F\left(\int_0^t \frac{r^{n-1}}{(1+r^2)^n} dr \right) = \int_0^t \frac{r^{n-2}}{(1+r^2)^{n-1}} dr, \quad t\ge 0.$$

\begin{Prop}
Let $K$ be a star body in $\mathbb S^n_+$. Then
$$\int_{S^{n-1}}  \mathrm{vol}(K\cap \xi^\perp)   d\xi\le 2^{n-1} |S^{n-1}||S^{n-2}|  F\left(   \frac{1}{|S^{n-1}| } \mathrm{vol}(K) \right) ,$$
with equality if and only if $K$ is a ball centered at the origin.

\end{Prop}

\proof It is not hard to check that $F$ is concave. Thus,
\begin{align*} \int_{S^{n-1}}  & \mathrm{vol}(K\cap \xi^\perp)   d\xi \\
& = 2^{n-1} \int_{S^{n-1}} \int_{S^{n-1}\cap \xi^\perp} \int_0^{\|\theta\|_K^{-1}} \frac{r^{n-2}}{(1+r^2)^{n-1}} dr \, d\theta\, d\xi \\
 & = 2^{n-1} |S^{n-2}| \int_{S^{n-1}}  \int_0^{\|\theta\|_K^{-1}} \frac{r^{n-2}}{(1+r^2)^{n-1}} dr \, d\theta \\
&= 2^{n-1} |S^{n-2}| \int_{S^{n-1}}   F\left(\int_0^{\|\theta\|_K^{-1}} \frac{r^{n-1}}{(1+r^2)^n} dr \right)  \, d\theta \\
& \le 2^{n-1} |S^{n-1}||S^{n-2}|  F\left(  \frac{1}{|S^{n-1}| } \int_{S^{n-1}} \int_0^{\|\theta\|_K^{-1}} \frac{r^{n-1}}{(1+r^2)^n} dr \, d\theta \right) \\
& = 2^{n-1} |S^{n-1}||S^{n-2}|  F\left(   \frac{1}{|S^{n-1}| } \mathrm{vol}(K) \right).
\end{align*}

\qed

Next we will prove a version of Busemann's intersection inequality in $\mathbb S^n_+$, which is, however, not optimal. %For centered balls the inequality is strict.

\begin{Prop}\label{Thm:not-opt}
Let $K$ be a star body in $\mathbb S^n_+$. Then
$$  \int_{S^{n-1}}   \mathrm{vol}(K\cap \xi^\perp)  ^n  d\xi < \frac{2^{n-1} n \kappa_{n-1}^n}{\kappa_n^{n-2}} \left[ \mathrm{vol}(K) \right]^{n-1}.$$

\end{Prop}

\proof
It was shown in \cite{DPP} that if $f$ is a non-negative bounded integrable function on $\mathbb R^n$, then
$$  \int_{S^{n-1}} \frac{\left( \int_{\xi^\perp} f(x) dx\right)^n}{\|f|_{\xi^\perp}\|_\infty}\,  d \xi \le \frac{n \kappa_{n-1}^n}{\kappa_n^{n-2}} \left( \int_{\mathbb R^n} f(x) dx \right)^{n-1}.$$

Now regard $K$ as a star body in $B^n$ and use the latter formula with  $$f(x) = \chi_K(x) \frac{2^{n-1}}{(1+|x|^2)^{n-1}}.$$
Since
$$\|f|_{\xi^\perp}\|_\infty = 2^{n-1},$$ for every $\xi\in S^{n-1}$, we get
$$\frac{1}{2^{n-1}} \int_{S^{n-1}}  \left[ \int_{K\cap \xi^\perp} \frac{2^{n-1}\, dx}{(1+|x|^2)^{n-1}} \right]^n d\xi  \le \frac{n \kappa_{n-1}^n}{\kappa_n^{n-2}} \left[ \int_{K}\frac{2^{n-1}\, dx}{(1+|x|^2)^{n-1}}\right]^{n-1},$$
that is
$$  \int_{S^{n-1}}  \left[\mathrm{vol}(K\cap \xi^\perp) \right]^n \,  d\xi  \le \frac{2^{n-1} n \kappa_{n-1}^n}{\kappa_n^{n-2}} \left[ \int_{K}\frac{2^{n-1}\, dx}{(1+|x|^2)^{n-1}} \right]^{n-1}.$$

Let us now estimate the right-hand side.
$$\int_{K}\frac{2^{n-1}}{(1+|x|^2)^{n-1}} dx \le \int_{K}\frac{2^{n}}{(1+|x|^2)^{n}} dx= \mathrm{vol}(K).$$
Finally,

$$  \int_{S^{n-1}}   \mathrm{vol}(K\cap \xi^\perp)  ^n  d\xi  \le \frac{2^{n-1} n \kappa_{n-1}^n}{\kappa_n^{n-2}} \left(  \mathrm{vol}(K) \right)^{n-1}.$$

\qed

\section{Busemann's intersection inequality in $\mathbb{S}^2_+$}

From the proof of  Theorem \ref{Thm:not-opt} it is clear that the
inequality obtained therein is not optimal. Our goal is to establish
an optimal inequality and we will start with the case of
$\mathbb{S}^2_+$. Here we will obtain sharp estimates for the minimum
and maximum of $\int_{S^1}  \mathrm{vol}(K\cap\xi^\perp)^2\, d\xi$ in
the class of origin-symmetric star or convex bodies (of a fixed area)
in $\mathbb{S}^2_+$. The case of $\mathbb H^2$ might suggest that
natural candidates for the maximizers are centered balls. However, this is
not true.
We will start with the following surprising result saying that centered balls in
$\mathbb{S}^2_+$ are in fact the  minimizers.

\begin{Thm}\label{th:min2d}
Let $K$ be an origin-symmetric star body in $\mathbb S^2_+$. Then
$$\int_{S^1}  \mathrm{vol}(K\cap\xi^\perp)^2\, d\xi\ge 8\pi\arccos^2\big(1-\mathrm{vol}(K)/2\pi\big),$$
with equality if and only if $K$ is a ball centered at the origin.
\end{Thm}

\proof Let $B$ be the centered ball in $\mathbb S^2_+$ of the same volume as $K$, and let $r\in (0,\pi/2)$ be the radius of $B$. Then
$$\mathrm{vol}(K)=\int_{S^1}\int_0^{\rho_K(u)}\sin t\,dt\,du=\int_{S^1}\left[1-\cos(\rho_K(u))\right] \,du$$
and
$$\mathrm{vol}(B)=\int_{S^1}\left[1-\cos(\rho_B(u))\right] \,du=2\pi(1-\cos r),$$
which yields $r=\arccos\big(1-\mathrm{vol}(K)/2\pi)$ and also
$$\int_{S^1}\left[\cos r-\cos(\rho_K(u))\right]\,du=0.$$
The desired inequality comes from the following inequality. For any $x,r\in(0,\pi/2]$ we have
\begin{equation}\label{eq:minineq}
x^2-r^2\ge \frac{2r}{\sin r}\big(\cos r-\cos x\big),
\end{equation}
with equality if and only if $x=r$. Indeed, the function
$$f(x)=x^2-r^2-\frac{2r}{\sin r}(\cos r-\cos x)$$ has the derivative $$f'(x)=2x-\frac{2r}{\sin r}\sin x.$$  The latter function is strictly convex in $(0,\pi/2]$ and satisfies
$f'(0)=f'(r)=0$. Thus $f'(x)<0$ when $0<x<r$ and $f'(x) > 0$ when $r<x<\pi/2$. Therefore,  $f$ has a strict minimum at $x=r$, i.e., $f(x)>f(r)=0$ for each $x\in(0,\pi/2]$ with $x\neq r$. Finally we have
\begin{align*}
\int_{S^1}\mathrm{vol}(K\cap\xi^\perp)^2d\xi &- \int_{S^1}\mathrm{vol}(B\cap\xi^\perp)^2 d\xi = 4\int_{S^1}\left[\rho_K^2(u)-r^2\right] \,du \\
&\ge \frac{8r}{\sin r}\int_{S^1}\left[\cos r-\cos(\rho_K(u))\right]\,du=0
\end{align*}
with the desired case of equality. The required inequality is now obtained by using $\mathrm{vol}(B\cap\xi^\perp)=2\arccos\big(1-\mathrm{vol}(K)/2\pi)$.

\qed

We will now turn to the maximizers.
\begin{Thm}
Let $K\subset \mathbb S^2_+$ be origin-symmetric and star-shaped. Then
$$\int_{S^1}  \mathrm{vol}(K\cap\xi^\perp)^2\, d\xi\le \pi^2\mathrm{vol}(K),$$
with equality if and only if $K$ is a spherical cone.
\end{Thm}

\proof
Let $C$ be an origin-symmetric cone   in $\mathbb S^2_+$ that has the same volume as $K$. Since $\rho_C$ takes the value $\frac\pi2$ on its support $A\subset S^1$, we get $\mathrm{vol}(C)=\int_A [1-\cos(\frac\pi2)]\,du=|A|$, and hence
$$\int_{S^1}\mathrm{vol}(C\cap\xi^\perp)^2d\xi =4\int_{S^1}\rho_C^2(\xi)d\xi=4(\pi/2)^2|A|=\pi^2\mathrm{vol}(K).$$
Thus,
\begin{align*}
&\int_{S^1}\mathrm{vol}(C\cap\xi^\perp)^2d\xi - \int_{S^1}\mathrm{vol}(K\cap\xi^\perp)^2 d\xi \\
&=\pi^2\int_{S^1}\left[1-\cos(\rho_K(u))\right] \,du -4\int_{S^1} \rho_K^2(u)  \,du\\
&=2\pi^2\int_{S^1} \left(\big[\sin(\rho_K(u)/2)\big]^2-\big[(\sqrt{2}/\pi)\rho_K(u)\big]^2\right) \,du.
\end{align*}
Finally, the inequality   $$\sin(x/2)\ge(\sqrt{2}/\pi)\, x ,\quad 0\le x\le\pi/2,$$ with its equality case,  $x=0$ or $x=\pi/2$, completes the proof.

\qed

In the next theorem we obtain an analogue of the previous result in
the class of origin-symmetric {\em convex} bodies  in $\mathbb{S}^2_+$. Here the maximizers are different; they are origin-symmetric lunes.

\begin{Thm}\label{th:2dstrip}
Let $K$ be an origin-symmetric convex body in $\mathbb S^2_+$. Then
$$\int_{S^1}  \mathrm{vol}(K\cap\xi^\perp)^2\, d\xi\le 16\int_0^{\pi/2}\arctan^2\Big(\frac{\tan(\mathrm{vol}(K)/4)}{\cos\theta}\Big)\,d\theta,$$
with equality if and only if $K$ is a lune.
\end{Thm}

To prove this theorem we need several auxiliary results. First, let us introduce some definitions.  If $u$ and $v$ are two vectors on $S^1$, then   $[u,v]$ stands for the segment of $S^{1}$ connecting $u$ and $v$, traced in the counterclockwise direction starting from $u$. By $\mathrm{cone}(u,v)$ we denote the spherical cone with base $[u,v]\subset S^1$. In particular, if $u=v$, then $\mathrm{cone}(u,u)$ stands for the geodesic ray from the origin in the direction of $u\in S^1$.

Let $x_0\in S^1$, and $K\subset\mathbb{S}^2_+$ be such that the interior of $K$ is connected and has the same volume as $K$.
For our purposes $K$ will be a connected component of the difference of two star bodies.
Define the function $f$ by
$$f(x)=\frac{\mathrm{vol}(K\cap\mathrm{cone}(x_0,x))}{\mathrm{vol}(K)}, \quad x\in D,$$
where $D\subset S^1$ is the closure of the set of all $x$'s such that the ray $\mathrm{cone}(x,x)$ intersects the interior of $K$.
Choosing $x_0\in S^1\setminus D$ if $S^1\setminus D\ne \emptyset$ or choosing any $x_0\in S^1$ otherwise, one can see that the function $f$ becomes injective. In this case the inverse function $f^{-1}:[0,1]\to S^1$ is called the {\em inverse angular area} function of $K$  measured counterclockwise from $x_0$. Note that $f$ and $f^{-1}$ do not depend on the choice of $x_0\in S^1\setminus D$ if $S^1\setminus D\ne \emptyset$. So in such situations we will just use the term the inverse angular area without specifying a direction $x_0$, but assuming that it lies outside of $D$.

\begin{Lem}\label{lem:icaf2d}
Let $K$, $\tilde K$ be star bodies in $\mathbb S^2_+$ such that $\tilde K\backslash K$ and $K\backslash\tilde K$ are connected and have equal positive volume. Let $\zeta^+$, $\zeta^-$ be the inverse angular area functions of $\tilde K\backslash K$, $K\backslash\tilde K$ correspondingly. Then
\begin{equation*}
\int_{S^1} \big[\rho_{\tilde K}^2 -\rho_K^2 \big](\theta)\, d\theta  > \mathrm{vol}(K\backslash\tilde K)\int_0^1 \Big[ F\big(\rho_K(\zeta^+(t))\big) - F\big(\rho_K(\zeta^-(t))\big)\Big] dt,
\end{equation*}
where $F(x)=2x/\sin x$. In particular, if  $$\rho_K(\zeta^+(t))\ge\rho_K(\zeta^-(t)),\quad t\in[0,1],$$ then $$\int_{S^1}\rho_{\tilde K}^2(\theta)  \, d\theta > \int_{S^1}\rho_K^2 (\theta) \, d\theta .$$
\end{Lem}

\proof Let both inverse angular area functions $\zeta^+$, $\zeta^-$ be measured counterclockwise from some point $x_0$ on $S^1$ identified with zero in $[0,2\pi)$. Note that
\begin{align*}
\mathrm{vol}(\tilde K\backslash K)=\int_0^{2\pi}\Big[(1-\cos\rho_{\tilde K})-(1-\cos\rho_K)\Big]_+\,d\theta =\int_0^{2\pi}g_+(\theta)\,d\theta
\end{align*}
and similarly $$\mathrm{vol}(K\backslash\tilde K)=\int_0^{2\pi}g_-(\theta)\,d\theta,$$ where  $g=\cos\rho_K-\cos\rho_{\tilde K}$ and $g_+=\max(g,0)$, $g_-=\max(-g,0)$. Furthermore, for each $t\in[0,1]$,
$$\int_0^{\zeta^+(t)}g_+(\theta)\,d\theta=t\,\mathrm{vol}(\tilde K\backslash K)$$
and
$$\int_0^{\zeta^-(t)}g_-(\theta)\,d\theta=t\,\mathrm{vol}(K\backslash\tilde K).$$
The functions $\zeta^+$, $\zeta^-$ are increasing
 and $\mathrm{vol}(K\backslash\tilde K)=\mathrm{vol}(\tilde K\backslash K)$,
so the above equalities can be written in a differential form:
\begin{equation}\label{eq:changevar}
g_+(\zeta^+)\,d\zeta^+=\mathrm{vol}(K\backslash\tilde K)\,dt\quad\text{and}\quad g_-(\zeta^+) \,d\zeta^-=\mathrm{vol}(K\backslash\tilde K)\,dt.
\end{equation}
Using   inequality  \eqref{eq:minineq}, we have
\begin{align*}
\int_0^{2\pi} \big[\rho_{\tilde K}^2(\theta)-\rho_K^2(\theta)\big]\,d\theta > \int_0^{2\pi} \frac{2\rho_K(\theta)}{\sin\rho_K(\theta)}\big[\cos\rho_K(\theta)-\cos \rho_{\tilde K}(\theta)\big]\,d\theta \\
=\int_{\mathrm{supp}(g_+)}F(\rho_K(\theta))g_+(\theta)\,d\theta - \int_{\mathrm{supp}(g_-)}F(\rho_K(\theta))g_-(\theta)\,d\theta
\end{align*}
where $F(x)=2x/\sin x$. The above inequality is strict because $K$ and $\tilde K$ are different. Also, the function $F$ is strictly increasing on $(0,\pi/2)$ because $F'(x)=2\cos x(\tan
x-x)/\sin^2x$ is positive for each $x\in(0,\pi/2)$.
Changing variables $\theta=\zeta^{\pm}(t)$ for each integral in the right-hand side of the equality above and using \eqref{eq:changevar}, we get
\begin{align*}
\int_{\mathrm{supp}(g_+)}F(\rho_K(\theta))g_+(\theta)\,d\theta &= \int_0^1F(\rho_K(\zeta^+(t)))g_+(\zeta^+(t))\,d\zeta^+(t) \\
&=\mathrm{vol}(K\backslash\tilde K)\int_0^1F(\rho_K(\zeta^+(t)))\,dt
\end{align*}
and
\begin{align*}
\int_{\mathrm{supp}(g_-)}F(\rho_K(\theta))g_-(\theta)\,d\theta &= \int_0^1F(\rho_K(\zeta^-(t)))g_-(\zeta^-(t))\,d\zeta^-(t) \\
&=\mathrm{vol}(K\backslash\tilde K)\int_0^1F(\rho_K(\zeta^-(t)))\,dt,
\end{align*}
which completes the proof.

\qed

%******************* Not edited below ************
Below when we say that we take a point on a line (or, more generally, on a curve) $\ell$ in the direction of a vector $u\in S^{1}$, it means that we consider the geodesic ray emanating from the origin in the direction of $u$ and take the point of intersection of this geodesic ray and the line $\ell$.

\begin{Lem}\label{lem:iaafstab}
Let $K,L\subset\mathbb S^2_+$ be convex bodies that lie outside of a centered ball of radius $r>0$ and such that $\mathrm{vol}(K\cap L)\ne 0$. Then the inverse angular area functions $\zeta_K$, $\zeta_L$  of $K$, $L$  satisfy
\begin{equation*}
\left\|\zeta_K-\zeta_L\right\|_{C[0,1]}^2 \le c\,\frac{\mathrm{vol}(K\Delta L)}{\mathrm{vol}(K\cap L)},
\end{equation*}
where $c=c(r)>0$ and $K\Delta L=(K\backslash L)\cup(L\backslash K)$.

\end{Lem}

\proof
If $$\frac{\mathrm{vol}(K\Delta L)}{\mathrm{vol}(K\cap L)}>\frac12,$$ then we are done since $\left\|\zeta_K-\zeta_L\right\|_{C[0,1]}^2\le \left( \|\zeta_K\|+\|\zeta_L\|\right)^2 = 4$ and we can take $c=8$.
Thus, it is enough to consider the case $\mathrm{vol}(K\Delta L)/\mathrm{vol}(K\cap L)\le 1/2$. We will  assume that $\zeta_K$, $\zeta_L$ are measured counterclockwise from some point $x_0\in S^1$, such that the geodesic ray in this direction does not intersect $K\cup L$. As before, we will identify $S^1$ with $[0,2\pi)$, and $x_0$ with $0\in [0,2\pi)$.

First, suppose that $\zeta_K(t)\ge\zeta_L(t)$ for a fixed $t\in[0,1]$. Then
$$\mathrm{vol}[K\cap\mathrm{cone}(0,\zeta_K(t))] = t\, \mathrm{vol}(K)$$
and
\begin{align*}
\mathrm{vol}[K\cap\mathrm{cone}(0,\zeta_L(t))] &\ge\mathrm{vol}[L\cap\mathrm{cone}(0,\zeta_L(t))] - \mathrm{vol}(L\backslash K) \\
&= t\,\mathrm{vol} (L) - \mathrm{vol} (L\backslash K),
\end{align*}
so
\begin{align}\label{eq:iaafstab_A}
&\mathrm{vol}[K\cap\mathrm{cone}(\zeta_L(t),\zeta_K(t))] \\
&= \mathrm{vol}[K\cap\mathrm{cone}(0,\zeta_K(t))] - \mathrm{vol}[K\cap\mathrm{cone}(0,\zeta_L(t))]\notag \\
&\le t\,\mathrm{vol} (K) - t\,\mathrm{vol} (L) + \mathrm{vol} (L\backslash K)\notag \\
&\le\mathrm{vol}(K\Delta L).
\end{align}
Using  \eqref{eq:iaafstab_A} and the assumption  $\mathrm{vol}(K\cap L)\ge 2\,\mathrm{vol}(K\Delta L)$, we get
\begin{align*}
\mathrm{vol}[K\setminus\mathrm{cone}(\zeta_L(t),\zeta_K(t))] &= \mathrm{vol}(K) - \mathrm{vol}[K\cap\mathrm{cone}(\zeta_L(t),\zeta_K(t))]\\
&\ge \mathrm{vol}(K\cap L) - \mathrm{vol}(K\Delta L) \ge \mathrm{vol}(K\Delta L)\\
&\ge \mathrm{vol}[K\cap\mathrm{cone}(\zeta_L(t),\zeta_K(t))].
\end{align*}

Let $p_0,p_1\in\partial K$ be in the direction of $\zeta_K(t)$ with $\mathrm{d}(p_0,o)\le\mathrm{d}(p_1,o)$ and $q_0,q_1\in\partial K$ be in the direction of $\zeta_L(t)$ with $\mathrm{d}(q_0,o)\le\mathrm{d}(q_1,o)$. Consider the region $S$ enclosed by two lines through $p_0,q_0$ and through $p_1,q_1$. Then convexity of $K$ implies
\begin{align*}
\mathrm{vol}[K\cap\mathrm{cone}(\zeta_L(t),\zeta_K(t))]&\ge\mathrm{vol}[S\cap\mathrm{cone}(\zeta_L(t),\zeta_K(t))],\\
\mathrm{vol}[K\backslash\mathrm{cone}(\zeta_L(t),\zeta_K(t))]&\le\mathrm{vol}[S\backslash\mathrm{cone}(\zeta_L(t),\zeta_K(t))],
\end{align*}
and hence
\begin{align*}
\mathrm{vol}(K) &= \mathrm{vol}[K\cap\mathrm{cone}(\zeta_L(t),\zeta_K(t))]+ \mathrm{vol}[K\backslash\mathrm{cone}(\zeta_L(t),\zeta_K(t))]\\
&\le 2\,\mathrm{vol}[K\backslash\mathrm{cone}(\zeta_L(t),\zeta_K(t))] \le 2\,\mathrm{vol}[S\backslash\mathrm{cone}(\zeta_L(t),\zeta_K(t))] \\
&\le 2\,\mathrm{vol}(S).
\end{align*}
Thus,
\begin{equation}\label{eq:iaafstab_B}
\frac{\mathrm{vol}[K\cap\mathrm{cone}(\zeta_L(t),\zeta_K(t))]}{\mathrm{vol}(K)} \ge \frac{\mathrm{vol}[S\cap\mathrm{cone}(\zeta_L(t),\zeta_K(t))]}{2\,\mathrm{vol}(S)}.
\end{equation}

We will now bound the right-hand side of \eqref{eq:iaafstab_B}.  By $\ell(\theta)$ we will denote the line passing through   the point $o_*=p_0q_0\cap p_1q_1$  and making an angle $\theta$ with the line $p_0q_0$.

\begin{figure}[h!]
\centering
\includegraphics[scale=.9,page=4]{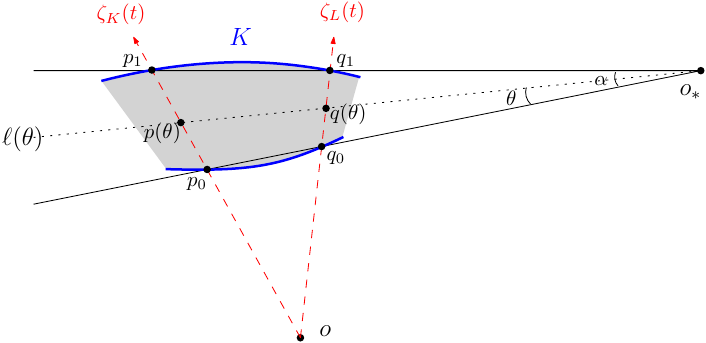}
\end{figure}
Let $p(\theta)$, $q(\theta)\in\ell(\theta)$ be the points in the direction of $\zeta_K(t)$, $\zeta_L(t)$ respectively.
If $\alpha$ is the angle between the lines through $p_0,q_0$ and through $p_1,q_1$, then $\mathrm{vol}(S)=2\alpha$ and
\begin{align*}
\mathrm{vol}[S\cap\mathrm{cone}(\zeta_L(t),\zeta_K(t))] &= \int_0^\alpha\Big|\cos\mathrm{d}(p(\theta),o_*)-\cos\mathrm{d}(q(\theta),o_*)\Big|d\theta \\
&\ge \frac2{\pi^2}\int_0^\alpha\mathrm{d}(p(\theta),q(\theta))^2 \,d\theta.
\end{align*}
The last inequality comes from $$\frac{|\cos x-\cos y|}2= \sin\Big(\frac{x+y}{2}\Big) \Big|\sin\Big(\frac{x-y}{2}\Big)\Big|\ge \sin^2\Big(\frac{x-y}{2}\Big)\ge\Big(\frac{x-y}{\pi}\Big)^2,$$ for $x,y\in[0,\pi/2]$.
Since, by the hypothesis of the lemma, $p(\theta), q(\theta)\in K$ are at least distance $r$  from the origin, we have
$$\mathrm{d}(p(\theta),q(\theta))=|\ell(\theta)\cap\mathrm{cone}(\zeta_L(t),\zeta_K(t))|\ge|\zeta_K(t)-\zeta_L(t)|\sin r.$$
Thus
\begin{align*} \frac{\mathrm{vol}[S\cap\mathrm{cone}(\zeta_L(t),\zeta_K(t))]}{2\,\mathrm{vol}(S)}&\ge\frac{\alpha(2/\pi^2)|\zeta_K(t)-\zeta_L(t)|^2\sin^2r}{4\alpha}\\
& =\frac{\sin^2r}{2\pi^2}|\zeta_K(t)-\zeta_L(t)|^2.
\end{align*}
Therefore, \eqref{eq:iaafstab_A} and \eqref{eq:iaafstab_B} imply $$|\zeta_K(t)-\zeta_L(t)|^2\le\frac{2\pi^2}{\sin^2r}\cdot\frac{\mathrm{vol}(K\Delta L)}{\mathrm{vol}(K)},$$
and similarly, in case of $\zeta_K(t)\le\zeta_L(t)$,
$$|\zeta_K(t)-\zeta_L(t)|^2\le\frac{2\pi^2}{\sin^2r}\cdot\frac{\mathrm{vol}(K\Delta L)}{\mathrm{vol}(L)},$$ which completes the proof.

\qed

%********** Not edited above*************

We will adopt the following notation. Each geodesic line $\ell$ in $\mathbb{S}^2_+$ not passing through the origin cuts $\mathbb{S}^2_+$ into two open sets with the common boundary $\ell$. These sets will be  denoted by $H^+_\ell$ and $H^-_\ell$,    where $H^-_\ell$  contains the origin and $H^+_\ell$ does not.

%For each geodesic line $\ell$ in $\mathbb{S}^2$ not passing through the origin we use the notations $H^+_\ell$ and $H^-_\ell$ for the open half-spaces in $\mathbb{S}^2$ with boundary $\ell$ %excluding the origin and including
%the origin, respectively.

\begin{Lem}\label{lem:icaf4lines}
Let $K\subset\mathbb S^2_+$ be a convex body and $\ell\subset\mathbb S^2_+$ be a line intersecting the interior of $K$, but not containing the origin. Consider a family of lines $\ell_\eps$ in $\mathbb S^2_+$ meeting $\ell$ at an angle $\eps$ (where  $\eps>0$ is small enough) such that
$$K_\eps=K\cap H^-_\ell\cap H^+_{\ell_\eps} \quad\left(\text{or } K_\eps=K\cap H^+_\ell\cap H^-_{\ell_\eps}\right)$$ is of positive measure and $\ell\cap\ell_\eps$ converges to a point $o_*\in\ell$ as $\eps\to 0$. Then the inverse angular area
function $\zeta_\eps$ of $K_\eps$ converges uniformly, as $\eps\to 0$, to a function  $\zeta:[0,1]\to S^1$ such that for each $t\in[0,1]$
\begin{equation*}
\cos\mathrm{d}(P_\ell\zeta(t),o_*)=(1-t)\cos\mathrm{d}(P_\ell\zeta(0),o_*) +t\cos\mathrm{d}(P_\ell\zeta(1),o_*).
\end{equation*}
Here, $P_\ell\xi$ denotes the point on $\ell$ in the  direction of $\xi\in S^1$.
\end{Lem}

\proof
Let $K_\eps=K\cap H^-_\ell\cap H^+_{\ell_\eps}$ be of positive measure for small $\eps>0$. The same argument works for the case of $K_\eps=K\cap H^+_\ell\cap H^-_{\ell_\eps}$.

For $0\le\theta\le\eps$, we denote by $\ell_\eps(\theta)$ the line through $o_\eps=\ell\cap\ell_\eps$ that intersects $\ell$ at an angle $\theta$ and also intersects $K_\eps$. (On the picture below the point $o_\eps$ is located outside of $K$, but in the proof there is no restriction on the position of $o_\eps$). For each $t\in[0,1]$ and $\theta\in[0,\eps]$,
consider the point $p_{t,\eps}(\theta)$ on the line $\ell_\eps(\theta)$ in the direction of $\zeta_\eps(t)$. Let $q_{t,\eps}(\theta)$ be the same point as $p_{t,\eps}(\theta)$ if $p_{t,\eps}(\theta)\in K$, and let $q_{t,\eps}(\theta)$ be the endpoint of
$K\cap \ell_\eps(\theta)$ closest to $p_{t,\eps}(\theta)$  if $p_{t,\eps}(\theta)\not\in K$. In case of $o_\eps\notin K$, the following figure shows where these points are located. If $o_\eps\in K$, all $p_{0,\eps}(\theta)$ and $q_{0,\eps}(\theta)$ indicate the common point $o_\eps$.

\begin{figure*}[h!]
\centering
\includegraphics[scale=.9,page=1]{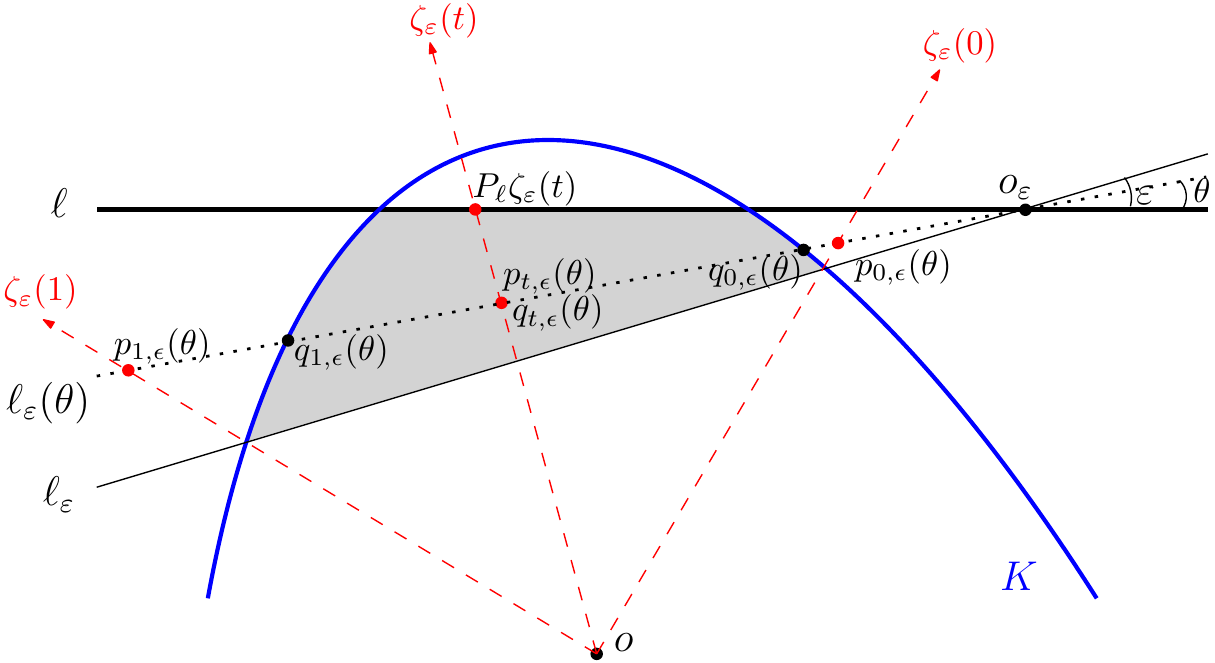}
\end{figure*}
Then the areas that appear in the definition of $\zeta_\eps$ can be computed by considering a polar coordinate system centered at $o_\eps$.  For all $t\in[0,1]$ we have
\begin{align*}
&\int_0^\eps\left(\big[1-\cos\mathrm{d}(q_{t,\eps}(\theta),o_\eps)\big] - \big[1-\cos\mathrm{d}(q_{0,\eps}(\theta),o_\eps)\big]\right)\,d\theta \\
&=t\int_0^\eps\left(\big[1-\cos\mathrm{d}(q_{1,\eps}(\theta),o_\eps)\big] - \big[1-\cos\mathrm{d}(q_{0,\eps}(\theta),o_\eps)\big]\right)\,d\theta.
\end{align*}
Rearranging, we have
\begin{align}\label{eq:2dstrip_lem}
&\frac1\eps\int_0^\eps\cos\mathrm{d}(q_{t,\eps}(\theta),o_\eps)\,d\theta \notag \\
&=\frac{1-t}\eps\int_0^\eps\cos\mathrm{d}(q_{0,\eps}(\theta),o_\eps)\,d\theta+\frac{t}\eps\int_0^\eps\cos\mathrm{d}(q_{1,\eps}(\theta),o_\eps)\,d\theta.
\end{align}
We first show that $$f_\eps(t)=\frac1\eps\int_0^\eps\cos\mathrm{d}\left(q_{t,\eps}(\theta),o_\eps\right)\,d\theta,\quad t\in[0,1],$$ converges uniformly (in $t$) as $\eps\to0$. Since $f_\eps(t)=(1-t)f_\eps(0)+tf_\eps(1)$, it suffices
to show that both $f_\eps(0)$ and $f_\eps(1)$ converge as $\eps\to0$.
Indeed, if $q_0$, $q_1$ are the endpoints of $K_\eps\cap\ell$ with $\mathrm{d}(q_0,o_*)<\mathrm{d}(q_1,o_*)$, then
\begin{align*}
\big|f_\eps(0)-\cos\mathrm{d}(q_0,o_*)\big| &\le\frac1\eps\int_0^\eps\big|\cos\mathrm{d}(q_{0,\eps}(\theta),o_\eps)-\cos\mathrm{d}(q_0,o_*)\big|\,d\theta \\
&\le \frac1\eps\int_0^\eps\big|\mathrm{d}(q_{0,\eps}(\theta),o_\eps) - \mathrm{d}(q_0,o_*)\big|\,d\theta \\
&\le \frac1\eps\int_0^\eps\left(\big|\mathrm{d}(q_{0,\eps}(\theta),o_\eps) - \mathrm{d}(q_0,o_\eps)\big|+\mathrm{d}(o_\eps,o_*)\right)  \,d\theta \\
&\le \mathrm{d}(q_{0,\eps}(\eps),q_0) + \mathrm{d}(o_\eps,o_*),
\end{align*}
and similarly $|f_\eps(1)-\cos\mathrm{d}(q_1,o_*)|\le\mathrm{d}(q_{1,\eps}(\eps),q_1) + \mathrm{d}(o_\eps,o_*)$. 
Since both $\mathrm{d}(q_{0,\eps}(\eps),q_0)$ and $\mathrm{d}(q_{1,\eps}(\eps),q_1)$ approach 0 as $\eps\to0$, it follows that $f_\eps(0)$ converges to $\cos\mathrm{d}(q_0,o_*)$, and $f_\eps(1)$ converges to $\cos\mathrm{d}(q_1,o_*)$.

Let $f$ be the limit of $f_\eps$ as $\eps\to0$. Define the function $\zeta:[0,1]\to S^1$ by
$$\cos\mathrm{d}(P_\ell\zeta(t),o_*)=f(t),\quad t\in[0,1].$$
Observing that $P_\ell\zeta_\eps(t)=p_{t,\eps}(0)$, we have
\begin{align*}
&\big|\cos\mathrm{d}(P_\ell\zeta_\eps(t),o_*)-\cos\mathrm{d}(P_\ell\zeta(t),o_*)\big| \\
&\le \big|\cos\mathrm{d}(P_\ell\zeta_\eps(t),o_*)-f_\eps(t)\big| + \|f_\eps-f\|_\infty \\
&=\frac1\eps\left|\int_0^\eps\big[\cos\mathrm{d}(p_{t,\eps}(0),o_*)-\cos\mathrm{d}(q_{t,\eps}(\theta),o_\eps)\big]\,d\theta\right| +\|f_\eps-f\|_\infty \\
&\le \frac1\eps\int_0^\eps\mathrm{d}\left(p_{t,\eps}(0),q_{t,\eps}(\theta)\right)\,d\theta + \mathrm{d}(o_\eps,o_*)+\|f_\eps-f\|_\infty.
\end{align*}
To show the uniform convergence of $\zeta_\eps$ to $\zeta$, it suffices to find a positive constant $c=c(K,\ell)$ such that
$$\mathrm{d}(p_{t,\eps}(0),q_{t,\eps}(\theta))\le c\eps,\quad0\le t\le1,\quad0\le\theta\le\eps.$$
Indeed, if $p_{t,\eps}(\theta)\in K$, then $q_{t,\eps}(\theta)=p_{t,\eps}(\theta)$ and hence $$\mathrm{d}(p_{t,\eps}(0),q_{t,\eps}(\theta))=\mathrm{d}(p_{t,\eps}(0),p_{t,\eps}(\theta))\le
\mathrm{d}(p_{t,\eps}(0),p_{t,\eps}(\eps)).$$
Otherwise, $q_{t,\eps}(\theta)$ is on the boundary of $K$, so we get
\begin{align*}
\mathrm{d}(p_{t,\eps}(0),q_{t,\eps}(\theta))&\le\mathrm{d}(p_{t,\eps}(0),\bar p)+\mathrm{d}(\bar p,q_{t,\eps}(\theta))\\
&\le\mathrm{d}(p_{t,\eps}(0),p_{t,\eps}(\eps))+\mathrm{d}(\ell\cap\ell_0,\ell_\eps\cap\ell_0),
\end{align*}
where $\ell_0$ is a supporting line to $K$ at $q_{t,\eps}(\theta)$ and $\bar p$ is the point  on the line $\ell_0$ in the direction of $\zeta_\eps(t)$.
\begin{figure}[h!]
\centering
\includegraphics[scale=1,page=2]{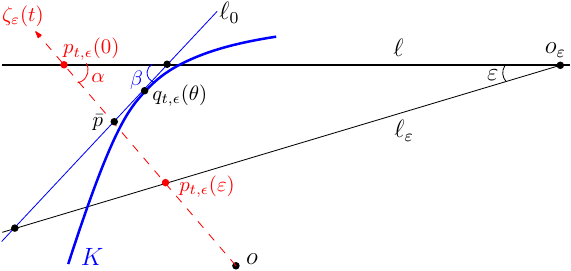}
\end{figure}

Then the law of sines for the triangle enclosed by $\ell$, $\ell_\eps$ and the direction of $\zeta_\eps(t)$ gives
$$\sin\mathrm{d}(p_{t,\eps}(0),p_{t,\eps}(\eps))=\sin\mathrm{d}(p_{t,\eps}(\eps),o_\eps)\cdot\frac{\sin\eps}{\sin\alpha}\le\frac{\sin\eps}{\sin\alpha},$$
and for the triangle enclosed by $\ell$, $\ell_\eps$ and $\ell_0$,
$$\sin\mathrm{d}(\ell\cap\ell_0,\ell_\eps\cap\ell_0)=\sin\mathrm{d}(\ell_\eps\cap\ell_0,o_\eps)\cdot\frac{\sin\eps}{\sin\beta}\le\frac{\sin\eps}{\sin\beta}.$$
Here $\alpha\in[0,\pi/2]$ is the angle between the line $\ell$ and the direction $\zeta_\eps(t)$, which is bounded below by a positive constant depending on $K\cap\ell$ only. Similarly, the angle $\beta$ between $\ell_0$ and
$\ell$ is bounded below by a constant depending only on $K$ and $\ell$. Thus both $\mathrm{d}(p_{t,\eps}(0),p_{t,\eps}(\eps))$ and $\mathrm{d}(\ell\cap\ell_0,\ell_\eps\cap\ell_0)$ are less than a constant multiple of $\eps$,
and so is $\mathrm{d}(p_{t,\eps}(0),q_{t,\eps}(\theta))$ in either case.

Therefore, $\zeta_\eps$ uniformly converges to $\zeta$, and we can get the desired equality for $\zeta$ by letting $\eps\to0$ in \eqref{eq:2dstrip_lem}.
\qed

\begin{Lem}\label{lem:midpoint2d}
Let $K$ be a convex body in $\mathbb S^2_+$ and $\ell$ a line in $\mathbb S^2_+$ intersecting the interior of $K$. For each small $\eps>0$, choose a line $\ell_\eps$ that makes an angle $\eps$ with $\ell$   and such that the regions of $K$  cut off by $\ell$ and $\ell_\eps$ have equal areas. Then the point $\ell\cap\ell_\eps$ converges to the midpoint of $K\cap \ell$ as $\eps\to0$.
\end{Lem}

\proof
We may assume that $\ell$ does not pass through the origin  and $K\cap H^+_\ell$ has the same volume as $K\cap H^+_{\ell_\eps}$. For $\theta\in[0,\eps]$, we denote by $\ell_\eps(\theta)$ the line that passes through $o_\eps=\ell\cap\ell_\eps$, makes an angle $\theta$ with  $\ell$,  and intersects $H^+_\ell\cap H^-_{\ell_\eps}$.

Let $p_\eps(\theta)$, $q_\eps(\theta)$ be the endpoints of $K\cap \ell_\eps(\theta)$ such that $p_\eps(\theta)\in H^+_\ell$ and $q_\eps(\theta)\in H^-_\ell$. In particular, if $\theta=0$, then $p_\eps(0)$ and $q_\eps(0)$ are the endpoints of the interval $K\cap \ell$. Write them as $p_0=p_\eps(0)$ and $q_0=q_\eps(0)$.
Then
\begin{align*}
0&=\frac1\eps\Big(\mathrm{vol}\big(K\cap H^+_\ell\big) - \mathrm{vol}\big(K\cap H^+_{\ell_\eps}\big)\Big) \\
&=\frac1\eps\int_0^\eps\left(\big[1-\cos\mathrm{d}(p_\eps(\theta),o_\eps)\big] - \big[1-\cos\mathrm{d}(q_\eps(\theta),o_\eps)\big]\right)\,d\theta \\
&=\frac1\eps\int_0^\eps\big[\cos\mathrm{d}(q_\eps(\theta),o_\eps)-\cos\mathrm{d}(p_\eps(\theta),o_\eps)\big]\,d\theta \\
&=\cos\mathrm{d}(q_0,o_\eps)-\cos\mathrm{d}(p_0,o_\eps)+g(\eps)-f(\eps),
\end{align*}
where
\begin{align*}
|f(\eps)|&=\frac1\eps\left|\int_0^\eps\big[\cos\mathrm{d}(p_\eps(\theta),o_\eps)-\cos\mathrm{d}(p_0,o_\eps)\big]\,d\theta\right| \\
&\le\frac1\eps\int_0^\eps\big|\mathrm{d}(p_\eps(\theta),o_\eps)-\mathrm{d}(p_0,o_\eps)\big|\,d\theta \le \mathrm{d}(p_\eps(\eps),p_0)
\end{align*}
and
\begin{align*}
|g(\eps)|&=\frac1\eps\left|\int_0^\eps\big[\cos\mathrm{d}(q_\eps(\theta),o_\eps)-\cos\mathrm{d}(q_0,o_\eps)\big]\,d\theta\right| \le \mathrm{d}(q_\eps(\eps),q_0).
\end{align*}
Since $\mathrm{d}(p_\eps(\eps),p_0) \to0$ and $\mathrm{d}(q_\eps(\eps),q_0) \to0$ as $\eps\to0$, we have that
$$\cos\mathrm{d}(p_0,o_\eps)-\cos\mathrm{d}(q_0,o_\eps)=g(\eps)-f(\eps)$$
converges to $0$ as $\eps\to0$. This implies that $o_\eps$ converges to the midpoint of $p_0$ and $q_0$ as $\eps\to0$.

\qed

\proof[{\bf Proof of Theorem~\ref{th:2dstrip}}] First we claim that there is a body that maximizes the integral $\int_{S^1} \rho_K^2(\theta)\, d\theta$ in the class of origin-symmetric convex bodies in $\mathbb S^2_+$ of a fixed volume. To see this, we can identify convex bodies in $\mathbb S^2_+$ with convex bodies in $\mathbb R^3$ that are obtained as the intersections of  the unit ball $B^3$ and convex cones in $\mathbb R^3$ with vertex at the origin. The claim  now follows from the Blaschke selection principle.

Suppose that there exists a body  $K$ that is not a lune and that maximizes $\int_{S^1} \rho_K^2(\theta)\, d\theta$ in the class of origin-symmetric convex bodies in $\mathbb S^2_+$ of a fixed volume. (We disregard the trivial case when the volume is equal to $\mathrm{vol}(\mathbb S^2_+ )$).  To get a contradiction we will construct an origin-symmetric convex
body $\tilde K\subset\mathbb S^2_+$ satisfying $\mathrm{vol}(K)=\mathrm{vol}(\tilde K)$ and $\int_{S^1} \rho_{\tilde K}^2(\theta)\, d\theta>\int_{S^1} \rho_K^2(\theta)\, d\theta$. We will consider the following cases according to the type of the boundary of $K$. By $[u,v]$ we denote the geodesic line segment between two points $u,v\in\mathbb S^2_+$. A line segment $[u,v]\subset\mathbb S^2_+$ is called an {\em edge} of $K$ if it coincides with the intersection of the boundary of $K$ and the line through $u$ and $v$.
\begin{enumerate}
\item $K$ has an arc on its boundary not containing any edges, i.e., the boundary of $K$ is strictly convex along the arc.
%\item there is a line that intersects the boundary of $K$ in a segment whose endpoints have different distances from the origin;
\item $K$ has an edge $[p_1^+,p_1^-]$ with $\mathrm{d}(p_1^+,o)>\mathrm{d}(p_1^-,o)$.
\item $K$ has two edges $[p_0,p_1]$ and $[q_0,q_1]$ such that $\mathrm{d}(p_0,p_1)<\mathrm{d}(q_0,q_1)$ and $\mathrm{d}(p_0,o)=\mathrm{d}(p_1,o)=\mathrm{d}(q_0,o)=\mathrm{d}(q_1,o)=\max_{x\in K}\mathrm{d}(x,o)$.
\item None of the above. In this case $K$ is reduced to a regular spherical polygon.
\end{enumerate}
It suffices to define such $\tilde K$ only in the half-space corresponding to $[0,\pi)\subset S^1$ and to assume that the arc and edges in the above cases lie on the half-space.

\medskip
\noindent{\bf Case (1)}: In this case we can choose two disjoint intervals $U,V\subset[0,\pi)$ such that
\begin{equation}\label{max-min}
\max_{u\in U}\rho_K(u)\le\min_{v\in V}\rho_K(v)
\end{equation}
 and the boundary of $K$ is strictly convex in the directions corresponding to  $U$ and $V$.
\begin{figure}[h!]
\centering
\includegraphics[scale=.7,page=3]{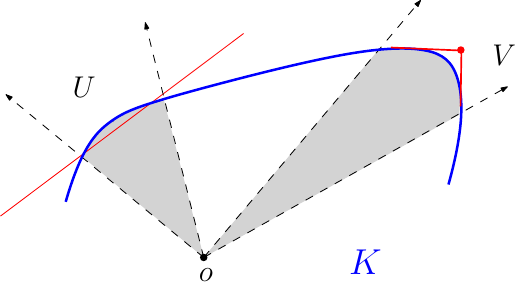}
\end{figure}
Strict convexity enables us to construct a convex body $\tilde K\subset\mathbb{S}^2_+$ such that $\mathrm{vol}(K)=\mathrm{vol}(\tilde K)$, and $\rho_K(x)\ge\rho_{\tilde K}(x)$ if $x\in U$, $\rho_K(x)\le\rho_{\tilde K}(x)$
if $x\in V$, and $\rho_K(x)=\rho_{\tilde K}(x)$ otherwise. Such a body $\tilde K$ can be obtained from $K$ by cutting off by a straight line a small piece of $K$ that lies in the region corresponding  to $U$, and then taking the convex hull of $K$ and a point outside of $K$ in the region corresponding to $V$. Moreover, we can ensure that  the volumes of $K$ and $\tilde K$ are the same. Strict convexity guarantees that the boundaries  of $K$ and $\tilde K$ can be different   only  in   the regions corresponding to $U$ and $V$. Using the notation of Lemma~\ref{lem:icaf2d}, we see that $\zeta^+$ and $\zeta^-$ take on their values only in $V$ and $U$ respectively. Therefore, by (\ref{max-min}) and  Lemma~\ref{lem:icaf2d} we have $\int_{S^1} \rho_{\tilde K}^2 \, d\theta>\int_{S^1} \rho_K^2(\theta)\, d\theta$.

\medskip
\noindent{\bf Case (2)}: Let $\ell$ be the line through $p_1^+$ and $p_1^-$. We will consider two possibilities: when $\ell$ is the only supporting line to the body $K$ at $p_1^+$, and when it is not. If
$\ell$ is the only supporting line to $K$ at $p_1^+$, we can use the same argument as in Case (1) to construct $\tilde K$. More precisely, we can cut off a small piece of $K$ in a neighborhood of $p_1^-$, and then take the convex hull with a point near $p_1^+$. The latter can be done as follows. Take a point $p$ on the line $\ell$ outside of $K$ and close to $p_1^+$. From this point $p$ draw the supporting line (different from $\ell$) to the body $K$. As $p$ approaches $p_1^+$, the supporting line approaches $\ell$ and the point of contact approaches $p_1^+$. Thus if we take such a point $p$ sufficiently close to $p_1^+$,  and take the convex hull of $p$ and the body $K$, then the resulting body will differ from $K$ only in a small neighborhood of $p_1^+$. The rest of the proof goes as in Case (1).

%  Note that $p_1^+$, $p_1^-$ are extreme points of $K$ because of maximality of the line segment. The uniqueness of a tangent line at $p_1^+$ implies that the tangent lines continuously change to $\ell$ as the
%tangency point approaches $p_1^+$, not through the line segment $p_1^+p_1^-$. Thus, if the convex hull is taken with a point on the line $\ell$ near $p_1^+$ but not in $K$, it changes $K$ only in a small neighborhood of
%$p_1^+$. Also, we can make the cut-off near $p_1^-$ as small as we want because $p_1^-$ is an extreme point; so $p_1^-$ is an exposed point itself or close to it. Therefore we can construct the desired $\tilde K$ in case
%that $K$ has only one tangent line at $p_1^+$.

Suppose now that there is another supporting line $\ell_0$ at $p_1^+$ other than $\ell$. For small $\eps>0$ choose a line $\ell_\eps$ satisfying the following three conditions: (i) $\ell_\eps$ separates $p_1^-$ from $p_1^+$ and
the origin; (ii) the angle between $\ell$ and $\ell_\eps$ is equal to $\eps$; (iii) the area of $K_\eps^+=H^-_{\ell_0}\cap  H^+_\ell\cap H^-_{\ell_\eps}$ is equal to that of $K_\eps^-=K\cap   H^+_{\ell_\eps}$.

\begin{figure}[h!]
\centering
\includegraphics[scale=.9,page=5]{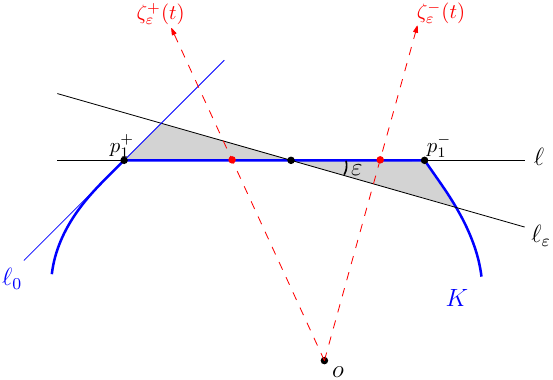}
\end{figure}
Then the body $\tilde K$ is obtained from $K$ by adding $K_\eps^+$ and by removing $K_\eps^-$. Let $\zeta_\eps^+$, $\zeta_\eps^-$ be the inverse angular area functions of $K_\eps^+$, $K_\eps^-$ measured from $\ell\cap\ell_\eps$ in opposite directions. Then, as $\eps\to0$, the point $\ell\cap\ell_\eps$ converges to the midpoint $p_0$ of $p_1^+$ and $p_1^-$ (Lemma~\ref{lem:midpoint2d}), and for
each choice of the sign $\pm$ the function $\zeta_\eps^\pm$ uniformly converges to  $\zeta^\pm$, by Lemma~\ref{lem:icaf4lines}, such that for each $t\in[0,1]$ the point $p_t^\pm\in\ell$ in the  direction of $\zeta^\pm(t)$ satisfies
\begin{align*}
\cos\mathrm{d}(p_t^\pm,p_0)&=(1-t)\cos\mathrm{d}(p_0^\pm,p_0)+t\cos\mathrm{d}(p_1^\pm,p_0) \\
&=1-t+t\cos\big(\mathrm{d}(p_1^+,p_1^-)/2\big).
\end{align*}
This implies that $\mathrm{d}(p_t^+,p_0)=\mathrm{d}(p_t^-,p_0)$ for all $t\in[0,1]$.
For each $t\in[0,1]$ and each of the signs $\pm$, the law of cosines for the triangle with vertices $o$, $p_0$ and $p_t^\pm$ yields
$$\cos\mathrm{d}(p_t^\pm,o)=\cos\mathrm{d}(p_0,o)\cos\mathrm{d}(p_t^\pm,p_0)\mp\sin\mathrm{d}(p_0,o)\sin\mathrm{d}(p_t^\pm,p_0)\cos\gamma,$$
where $\gamma$ is the angle between the lines through $o,p_0$ and through $p_0,p_1^+$. Note that $\gamma\ne \pi/2$. Since $\mathrm{d}(p_t^+,p_0)=\mathrm{d}(p_t^-,p_0)$,  the above equations imply
$$\cos\mathrm{d}(p_t^-,o)-\cos\mathrm{d}(p_t^+,o)=2\sin\mathrm{d}(p_0,o)\sin\mathrm{d}(p_t^+,p_0)\cos\gamma,$$
and therefore
\begin{align*}
\frac{\cos\mathrm{d}(p_t^-,o)-\cos\mathrm{d}(p_t^+,o)}{\cos\mathrm{d}(p_1^-,o)-\cos\mathrm{d}(p_1^+,o)}&=\frac{2\sin\mathrm{d}(p_0,o)\sin\mathrm{d}(p_t^+,p_0)\cos\gamma}{2\sin\mathrm{d}(p_0,o)\sin\mathrm{d}(p_1^+,p_0)\cos\gamma}
\\
&=\frac{\sin\mathrm{d}(p_t^+,p_0)}{\sin\mathrm{d}(p_1^+,p_0)}.
\end{align*}
Thus $$\mathrm{d}(p_t^+,o)-\mathrm{d}(p_t^-,o)\ge \cos\mathrm{d}(p_t^-,o)-\cos\mathrm{d}(p_t^+,o)=c\sin\mathrm{d}(p_t^+,p_0),$$
where  $c=\big({\cos\mathrm{d}(p_1^-,o)-\cos\mathrm{d}(p_1^+,o)}\big)/\sin\mathrm{d}(p_1^+,p_0)$ is a positive constant depending only on $K\cap\ell$.
Since $\rho_K(\zeta^+(t))=\mathrm{d}(p_t^+,o)$ and $\rho_K(\zeta^-(t))=\mathrm{d}(p_t^-,o)$, we have the strict inequality $$\rho_K(\zeta^+(t))>\rho_K(\zeta^-(t)),\quad t\in(0,1].$$
If $F$ is a strictly increasing function (as in Lemma~\ref{lem:icaf2d}), then $$\int_0^1F\big(\rho_K(\zeta^+(t))\big)\,dt>\int_0^1F\big(\rho_K(\zeta^-(t))\big)\,dt.$$ Since $\zeta_\eps\to\zeta$ uniformly as $\eps\to0$,
the same strict inequality holds for $\zeta_\eps$ instead of $\zeta$   for $\eps>0$  small enough. For such small $\eps>0$, apply Lemma~\ref{lem:icaf2d} to complete the proof in this case.

\medskip

%********************* Not edited below ***************************

\noindent{\bf Case (3)}:
Let $\ell^+$ be the line through $p_0$, $p_1$, and $\ell^+_\eps$ be the line that is parallel to $\ell^+$ but distance $\eps$ farther away from the origin. Similarly, $\ell^-$ stands for the line through $q_0$, $q_1$, and $\ell^-_\delta$   for the line that is parallel to $\ell^-$ but distance $\delta$ closer to the origin.

Let $B$ be the centered ball containing points $p_0,p_1,q_0,q_1$ on its boundary. Consider the subregions $B_\eps^+$ and $B_\delta^-$ of $B$ enclosed by $\ell^+$, $\ell^+_\eps$ and by $\ell^-$, $\ell^-_\delta$. Then $\tilde K$ is obtained from $K$ by adding the region $B_\eps^+$ and then removing the region outside the line $\ell_\delta^-$. It means that $\tilde K\backslash K=B_\eps^+\backslash B_\delta^-$ and $K\backslash \tilde K$ is the region enclosed by $\ell^-$, $\ell^-_\delta$, $K$. Here, $\delta=\delta(\eps)$ is chosen so that $\mathrm{vol}(\tilde K\backslash K)=\mathrm{vol}(K\backslash\tilde K)$.
\begin{figure}[h!]
\centering
\includegraphics[scale=0.8,page=6]{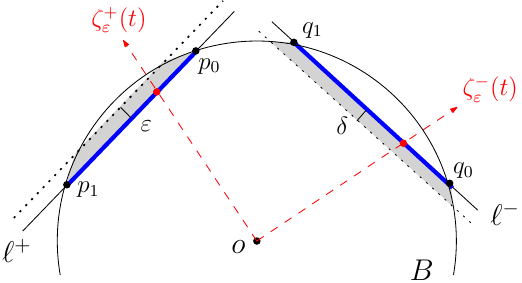}
\end{figure}

Let $\xi_\eps^+$, $\zeta_\eps^+$, $\zeta_\eps^-$ be the inverse angular area functions of $B_\eps^+$, $\tilde K\backslash K$, $K\backslash\tilde K$, respectively. If $\{p_0,p_1\}\cap\{q_0,q_1\}=\emptyset$, then $\xi_\eps^+=\zeta_\eps^+$ for small $\eps>0$. Otherwise, it follows from Lemma~\ref{lem:iaafstab} that $\|\xi_\eps^+-\zeta_\eps^+\|^2$ is less than a constant multiple of
$$\frac{\mathrm{vol}(B_\eps^+\Delta(\tilde K\backslash K))}{\mathrm{vol}(B_\eps^+\cap(\tilde K\backslash K))} \le \frac{\mathrm{vol}(B_\eps^+\cap B_\delta^-)}{\frac12 \mathrm{vol}(B_\eps^+)}\le\frac{c_1\eps\delta}{c_2\eps}=c_3\delta,$$
where the constants $c_1,c_2,c_3$ only depend on $\ell^+$, $\ell^-$, and $B$.

Applying Lemma~\ref{lem:icaf4lines} to the body $B$ and the line $\ell^+$ we have the following: as $\eps\to0$ the function $\xi_\eps^+$ (and thus $\zeta_\eps^+$ as well) uniformly converges to  $\zeta^+$ such that for each $t\in[0,1]$ the point $p_t\in\ell^+$ in the direction of $\zeta^+(t)$ satisfies
$$\cos\mathrm{d}(p_t,p_*)=(1-t)\cos\mathrm{d}(p_0,p_*)+t\cos\mathrm{d}(p_1,p_*)$$ where $p_*\in\ell^+$ lies on the boundary of $\mathbb S^2_+$.
Since
$$\big|\mathrm{d}(p_t,p_*)-\mathrm{d}(p_{1/2},p_*)\big|=\mathrm{d}(p_t,p_{1/2})$$
and
$$\mathrm{d}(p_{1/2},p_*)=\frac{\mathrm{d}(p_0,p_*)+\mathrm{d}(p_1,p_*)}2=\frac\pi2,$$
we have
\begin{align}\label{eq:sin_case3}
\sin\mathrm{d}(p_t,p_{1/2})&= \sin\big|\mathrm{d}(p_t,p_*)-\pi/2\big|=|\cos\mathrm{d}(p_t,p_*)| \notag\\
&=\big|(1-t)\cos\mathrm{d}(p_0,p_*)+t\cos\mathrm{d}(p_1,p_*)\big| \notag\\
&=|(1-2t)\cos\mathrm{d}(p_0,p_*)| =|1-2t|\sin\mathrm{d}(p_0,p_{1/2}).
\end{align}
The (spherical) Pythagorean theorem for the right triangle with vertices $o,p_{1/2},p_t$ yields
$$\cos\mathrm{d}(p_t,o)=\cos\mathrm{d}(p_{1/2},o)\cos\mathrm{d}(p_t,p_{1/2}),$$
and hence
$$\frac{\cos\mathrm{d}(p_t,o)}{\cos\mathrm{d}(p_0,o)}=\frac{\cos\mathrm{d}(p_t,p_{1/2})}{\cos\mathrm{d}(p_0,p_{1/2})}=\sqrt{\frac{1-\sin^2\mathrm{d}(p_t,p_{1/2})}{\cos^2\mathrm{d}(p_0,p_{1/2})}}.$$
Substituting  \eqref{eq:sin_case3} into the above equation, we get
\begin{equation*}\label{eq:pt_case3}
\cos\mathrm{d}(p_t,o) =\cos\mathrm{d}(p_0,o)\sqrt{1+4t(1-t)\tan^2\mathrm{d}(p_0,p_{1/2})}.
\end{equation*}
Similarly, the function $\zeta_\eps^-$ uniformly converges to $\zeta^-$ such that for each $t$ the point $q_t\in\ell^-$ in the direction of $\zeta^-(t)$ satisfies
\begin{equation*}\label{eq:qt_case3}
\cos\mathrm{d}(q_t,o)=\cos\mathrm{d}(q_0,o)\sqrt{1+4t(1-t)\tan^2\mathrm{d}(q_0,q_{1/2})}.
\end{equation*}
Since $\mathrm{d}(p_0,o)=\mathrm{d}(q_0,o)$ and $\mathrm{d}(p_0,p_{1/2})<\mathrm{d}(q_0,q_{1/2})$,
the above two equations imply
\begin{equation*}
\rho_K(\zeta^+(t))=\mathrm{d}(p_t,o)>\mathrm{d}(q_t,o)=\rho_K(\zeta^-(t)),\quad t\in(0,1).
\end{equation*}
As in Case (3), use the uniform convergence of $\zeta_\eps$ to $\zeta$ to get a strict inequality needed for Lemma~\ref{lem:icaf2d}.

\medskip
\noindent{\bf Case (4)}:
Suppose that an origin-symmetric convex body $K$ does not contain any type of arcs on its boundary corresponding to Cases (1) and (2).
As we will see below this implies that the boundary of $K$ consists of (possibly infinitely many) edges whose endpoints lie on the same centered circle, as well as their limiting points.
Consider the set $\mathcal I$ of all edges of $K$. Because of Case (2), we can assume that for any $I\in\mathcal I$ both endpoints of $I$ are   the same distance from the origin, which we will denote by $r(I)$. Let $R=\max_{u\in S^1}\rho_K(u)$. We claim that there is no edge $I$ with $r(I)<R$. To reach a contradiction, assume that there exists an edge $I_0$ with $r(I_0)<R$. We will show that  the sets $\mathcal I_{\epsilon} = \{I\in \mathcal I: R-\eps<r(I)<R\}$ are not empty (and therefore have infinitely many edges) for every $\epsilon >0$.
Indeed, suppose that $\mathcal I_{\epsilon}$ is empty for some $\epsilon>0$. Consider the following closed subsets of $\partial K$:  $A_1= \{x\in \partial K: d(x,o) = R\}$ and $A_2 = \{x\in \partial K: d(x,o) = R-\epsilon \}$.
Let $\gamma\subset \partial K$ be an arc of minimal length connecting $A_1$ and $A_2$. This arc must contain an edge, since we ruled out Case (1). But this edge cannot have $r(I)<R$, since $\mathcal I_{\epsilon}$ is empty  and for any interior point $x$ of $\gamma$ we have $d(x,o)>R-\epsilon$. Similarly, $\gamma$ cannot contain an edge with $r(I)=R$  since for any interior point $x$ of $\gamma$ we have $d(x,o)<R$. Contradiction.   Thus the number of edges  in $\mathcal I_{\epsilon}$  for  any $\eps>0$ should be   infinite. But this is also impossible.    Just fix $I\in\mathcal I$ with $r(I)<R$ and choose $J\in\mathcal I$ with $r(I)<r(J)<R$, whose length is small enough to ensure $\mathrm{d}(x,0)>r(I)$ for all $x\in J$. Use the same argument as in Case (1) for $I$, $J$ instead of the arcs corresponding to $U$, $V$; that is, to get $\tilde K$ cut off a small piece of $K$ around $I$ by a parallel line  and take a convex hull with a point close to $J$.
Thus we conclude that $r(I)=R$ for all $I\in\mathcal I$.
Finally, ruling out the arcs corresponding to Case (3), we see that   $K$ is reduced to a regular spherical polygon with $2m$ edges for some $m\ge2$.

Let $P$ be the body in $\mathbb R^2$ with radial function $\rho(x)=\tan\rho_K(x)$. Then $P\subset\mathbb R^2$ is a (Euclidean) regular $2m$-polygon centered at the origin. If $P$ has an inscribed circle of radius $r$, its radial function $\rho$ can be written as
$$\rho(x)=r\sec\left(x-\frac{\pi}{2m}\right),\quad x\in[0,\pi/m]$$
and   extended with period $\frac\pi{m}$ to other $x\in[0,2\pi]$.
Define the origin-symmetric body $\bar P\subset\mathbb R^2$ by its radial function $\bar\rho$ as follows.
\begin{equation*}
\bar\rho(x)=
\begin{cases}
\rho(\frac{x}2),&0\le x\le\frac{2\pi}m ,\\
\rho(x),&\frac{2\pi}m\le x\le\pi, \\
\bar\rho(x-\pi),&\pi\le x\le2\pi.
\end{cases}
\end{equation*}
Then we can use the formula of the curvature of a planar polar curve to show that  $\bar P$ is still convex and moreover its boundary contains a  strictly convex arc. Indeed, the signed curvature of $\partial\bar P$ at $x\in(0,\frac{2\pi}m)$ is
\begin{align*}
\bar\kappa(x)&= \frac{\bar\rho^2(x)+2(\bar\rho'(x))^2-\bar\rho(x) \bar\rho''(x)}{(\bar\rho^2(x)+(\bar\rho'(x))^2)^{3/2}}= \frac{\rho(\frac{x}2)^2+\frac12\rho'(\frac{x}2)^2-\frac14\rho(\frac{x}2)\rho''(\frac{x}2)}{(\bar\rho^2(x)+(\bar\rho'(x))^2)^{3/2}} \\
&=\frac34\cdot\frac{\bar\rho^2(x)}{(\bar\rho^2(x)+(\bar\rho'(x))^2)^{3/2}}+\frac14\cdot \frac{\rho(\frac{x}2)^2+2\rho'(\frac{x}2)^2-\rho(\frac{x}2)\rho''(\frac{x}2)}{(\bar\rho^2(x)+(\bar\rho'(x))^2)^{3/2}}
\end{align*}
Note that $\rho(\frac{x}2)^2+2\rho'(\frac{x}2)^2-\rho(\frac{x}2)\rho''(\frac{x}2)=0$ because the curvature of $\partial P$ at $x/2\in(0,\pi/m)$ is equal to 0. Thus,
\begin{equation*}
\bar\kappa(x)= \frac{3\bar\rho^2(x)}{4(\bar\rho^2(x)+(\bar\rho'(x))^2)^{3/2}}>0,\quad x\in(0,2\pi/m).
\end{equation*}
This means that the arc of $\bar P$ corresponding to $0<x<2\pi/m$ is strictly convex. To prove convexity of the entire body  $\bar P$ we need  to check what happens at the endpoints of this arc. It suffices to find a supporting line of $\bar P$ at each of these endpoints that does not intersect the interior of $\bar P$. This is true since   $\bar P$ is inscribed in the circle of radius
$r\sec\left( \frac{\pi}{2m}\right)$ and those endpoints lie on this circle.

Thus the body $\tilde K$ in $\mathbb{S}^2_+$ defined by $\tan\rho_{\tilde K}(x)=\bar\rho(x)$ is an origin-symmetric convex body such that
\begin{align*}
\mathrm{vol}(\tilde K)&=2\int_0^{\frac{2\pi}m}\big[1-\cos\rho_K(x/2)\big]\,dx +2(m-2)\int_0^{\frac\pi{m}}\big[1-\cos\rho_K(x)\big]\,dx \\
&=2m\int_0^{\frac\pi{m}}\big[1-\cos\rho_K(x)\big]\,dx = \mathrm{vol}(K)
\end{align*}
and $\int_0^{2\pi} \rho_{\tilde K}^2(x)\,dx=\int_0^{2\pi} \rho_K^2(x)\,dx$ similarly. Moreover, due to Case (1), $\tilde K$ cannot be a maximizer because it contains a  strictly convex arc. Thus $K$ is not a maximizer.

To finish the proof it remains to show that
$$\int_{S^1}  \mathrm{vol}(K\cap\xi^\perp)^2\, d\xi = 16\int_0^{\pi/2}\arctan^2\Big(\frac{\tan(\mathrm{vol}(K)/4)}{\cos\theta}\Big)\,d\theta,$$
if $K$ is an origin-symmetric lune.
This is true because the volume and the radial function of an origin-symmetric lune $K$ can be written as $$\mathrm{vol}(K) = 4w,\quad\tan\rho_K(\cdot) =\frac{\tan w}{|\langle\cdot,u\rangle|}$$
where $w\in(0,\frac\pi2)$ and $u\in S^1$ satisfy $w=\rho_K(u)=\min_{x\in K}\rho_K(x)$.

\qed

%*********** Edit start, Apr 15 *************

\section{Busemann's intersection inequality in $\mathbb{S}^n_+$,
  $n\ge3$}

In this section we study the maximizing and minimizing problems for the quantity
\begin{equation*}
\int_{S^{n-1}}\mathrm{vol}(K\cap \xi^\perp)^n d\xi
\end{equation*}
in the class of star bodies of a fixed volume in the spherical space
$\mathbb S^n_+$ when the dimension $n$ is $3$ or higher.

 Recall that in the hyperbolic space of any dimension  the centered balls are the unique maximizers of this
 quantity. On the other hand,   the centered
 balls in $\mathbb S^2_+$ are the unique minimizers. Thus it is quite
 surprising that the balls in  $\mathbb S^n_+$, $n\ge 3$, are neither
 minimizers nor maximizers (even in the class of origin-symmetric
 convex bodies).

\begin{Thm} The centered ball of any radius is neither a minimizer
  nor a maximizer of the integral $\int_{S^{n-1}} \mathrm{vol}(K\cap\xi^\perp)^n\,
  d\xi$   in the
  class of origin-symmetric convex bodies $K\subset\mathbb S^n_+$, $n\ge3$, of
  a fixed volume.
\end{Thm}

\proof
We will  show that  centered balls cannot even be local minimizers or
local maximizers of the integral in question.

Let $B\subset\mathbb S^n_+$ be the centered ball of radius
$r\in(0,\pi/2)$  and $K\subset\mathbb S^n_+$ be an origin-symmetric
convex body of the same volume as $B$. Let $f=\rho_K- r$ and
denote $\eps=\|f\|_{L^\infty(S^{n\!-\!1})}$ and
$\delta=\|f\|_{L^2(S^{n\!-\!1})}$. We will assume  $\eps$, and
therefore $\delta$,  to be sufficiently small. We have
$$\mathrm{vol}(K)=\int_{S^{n-1}}\int_0^{\rho_K(u)}\sin^{n-1}
t\,dt\,du=\int_{S^{n-1}}\phi_n(\rho_K(u)) \,du$$ and
$\mathrm{vol}(B)=\int_{S^{n-1}}\phi_n(r) \,du$, where
$\phi_n(x)=\int_0^x\sin^{n-1}t\,dt$. Using  the Taylor
expansion of $\phi_n$ at $x=r$, we get
$$\phi_n(x)-\phi_n(r)=\sin^{n-1}r\left[(x-r)+\frac{n-1}{2\tan r}(x-r)^2\right] + O(|x-r|^3).$$
Putting $\rho_K$ instead of $x$ in the above expansion and integrating both sides yield
\begin{align*}
0&=\mathrm{vol}(K)-\mathrm{vol}(B)=\int_{S^{n-1}}\big[\phi_n(\rho_K(u))-\phi_n(r)\big]\,du \\
&= \sin^{n-1}r\left(\int_{S^{n-1}}\big[f(u)+\frac{n-1}{2\tan r}f(u)^2\big]\,du\right) + O\left(\int_{S^{n-1}}|f(u)|^3du\right).
\end{align*}
Let $c_0=\frac{n-1}{2\tan r}$. Since
$$\int_{S^{n-1}}|f(u)|^3\, du \le \eps\int_{S^{n-1}}|f(u)|^2du= \eps \delta^2 ,$$ the above equality implies
\begin{align}\label{eq:locmaxeqvol}
\int_{S^{n-1}}f(u)\, du =-c_0\int_{S^{n-1}}f^2(u)\, du +O(\eps\delta^2) =-c_0\delta^2+O(\eps\delta^2).
\end{align}
We now compute the volumes of the sections of $K$. For each $\xi\in
S^{n-1}$ we have
 $$\mathrm{vol}(K\cap\xi^\perp)=\int_{S^{n-1}\cap\xi^\perp}\phi_{n-1}(\rho_K(u)) \,du=R(\phi_{n-1}\circ\rho_K)(\xi).$$
Expanding $\phi_{n-1}$ at $r$, we get
\begin{align*}
\mathrm{vol}(K\cap\xi^\perp) & =R\left[\phi_{n-1}(r)+\sin^{n-2}r\Big(f+\frac{n-2}{2\tan r}f^2\Big) + O(|f|^3)\right](\xi)\\
&= \mathrm{vol}(B\cap\xi^\perp)+ c_1R(f+c_2f^2)(\xi)+ O\big(R|f|^3(\xi)\big),
\end{align*}
where $c_1=\sin^{n-2}r$ and $c_2=\frac{n-2}{2\tan r}$.

Raising  both sides to the $n$th power, integrating over $S^{n-1}$, and using \eqref{eq:locmaxeqvol}, we have
\begin{align}
\label{ball3d}
&\int_{S^{n-1}}\mathrm{vol}(K\cap\xi^\perp)^n\,d\xi-
  \int_{S^{n-1}}\mathrm{vol}(B\cap\xi^\perp)^n\,d\xi \notag \\
 &=\int_{S^{n-1}} \Big(c_3\big[c_1R(f +c_2f^2)\big]+c_4\big[c_1R(f+c_2f^2)\big]^2\Big)(\xi)\, d\xi + O(\eps\delta^2)\notag \\
&=c_1c_3|S^{n-2}|\int_{S^{n-1}} (f+c_2f^2)(\xi)\,
  d\xi+ c_1^2c_4\int_{S^{n-1}} (Rf)^2(\xi) \, d\xi + O(\eps\delta^2)\notag \\
&=-c_1c_3(c_0-c_2)|S^{n-2}|\int_{S^{n-1}} f^2(\xi)\,
  d\xi+ c_1^2c_4\int_{S^{n-1}} (Rf)^2(\xi) \, d\xi + O(\eps\delta^2)\notag \\
&=-c_1^2c_4\int_{S^{n-1}} \Big[c_5f^2(\xi)- (Rf)^2(\xi)\Big]\, d\xi + O(\eps\delta^2),
\end{align}
where $c_3=n\mathrm{vol}(B\cap\xi^\perp)^{n-1}$, $c_4=\frac{n(n-1)}{2}\mathrm{vol}(B\cap\xi^\perp)^{n-2}$, and
\begin{align*}
c_5=\frac{c_1c_3(c_0-c_2)|S^{n-2}|}{c_1^2c_4}=\frac{|S^{n-2}|\,\mathrm{vol}(B\cap\xi^\perp)}{(n-1)\tan r\sin^{n-2}r}<\frac{|S^{n-2}|^2}{(n-1)^2}.
\end{align*}
The above inequality comes from $$\frac{\mathrm{vol}(B\cap\xi^\perp)}{|S^{n-2}|}=\int_0^r\sin^{n-2}t\,dt<\frac1{\cos
  r}\int_0^r\cos t\sin^{n-2}t\,dt=\frac{\tan r\sin^{n-2}r}{n-1}.$$
The remainder term in \eqref{ball3d} is obtained from the following  relations:
\begin{align*}
\int_{S^{n-1}} R|f|^3(\xi)\,d\xi =|S^{n-2}|\int_{S^{n-1}}|f|^3(\xi)\, d\xi = O(\eps \delta^2),
\end{align*}
\begin{align*}
\left|\int_{S^{n-1}}Rf(\xi)  R f ^2(\xi)\,d\xi\right| & \le \eps \int_{S^{n-1}}(R|f|(\xi))^2  \,d\xi \\
&  \le \eps |S^{n-2}|^2 \int_{S^{n-1}}f^2(\xi)\, d\xi = O(\eps \delta^2),
\end{align*}
and
\begin{align*}
&\int_{S^{n-1}}|R(f+c_2f^2)|^3 (\xi)\,d\xi \le(1+c_2\eps)^3 \int_{S^{n-1}} R|f|^3 (\xi)\,d\xi\\
&\le\eps(1+c_2\eps)^3|S^{n-2}| \int_{S^{n-1}}(R|f|)^2 (\xi)\,d\xi \\
&\le\eps(1+c_2\eps)^3|S^{n-2}|^3 \int_{S^{n-1}}f^2 (\xi)\,d\xi=O(\eps\delta^2),
\end{align*}
where we used  \eqref{L2bound} and \eqref{self-adj}.

To finish the proof, we provide  examples of two origin-symmetric
convex bodies for which the expression in
\eqref{ball3d} has opposite signs.
%Let $H_k$ be a spherical harmonic of an even degree $k$ with $\|H_k\|_{L^2(S^{n-1})}=1$ where $k$ will be specified later. Consider the body $K\subset\mathbb S^n$ with the radial function $$\rho_K(\xi)= r + \alpha + \beta H_k(\xi), \qquad \xi\in S^{n-1},$$ where $\alpha$, $\beta$ are chosen to satisfy $\alpha^2+\beta^2=\delta^2$ and $\mathrm{vol}(K)=\mathrm{vol}(B)$. Then
 Let $$f=\alpha  +\beta H_k,$$ where $H_k$ is a spherical harmonic of an
 even degree $k>0$ with $\|H_k\|_{L^2(S^{n-1})}=1$, and $\alpha$, $\beta$
 are  the numbers satisfying the conditions $\delta=\|f\|_{L^2(S^{n-1})}$ and $\mathrm{vol}(K)=\mathrm{vol}(B)$. The first condition implies $|S^{n-1}| \alpha^2+\beta^2=\delta^2$, and the second one gives
$$\alpha = \frac1{|S^{n-1}|}\int_{S^{n-1}} f(\xi)\, d\xi =-\frac{c_0\delta^2}{|S^{n-1}|} + O(\eps\delta^2),$$
which follows from  condition \eqref{eq:locmaxeqvol}.
Note that $\alpha\to 0$ and $\beta\to 0$ as $\delta\to0$. Choosing $\delta$ sufficiently
small, we can guarantee that the body $K$ is convex.

Now compute the sign of (\ref{ball3d}) for such a body $K$.
Since $RH_k= \lambda_k H_k$, $\int_{S^{n-1}}H_k=0$ and $\alpha=O(\delta^2)$, we  have
\begin{align*}
&\int_{S^{n-1}} \Big[c_5f^2(\xi)- (Rf)^2(\xi)\Big]\,d\xi  \\
&=\int_{S^{n-1}}\Big[c_5\big(\alpha^2 + \beta^2H_k^2(\xi)\big) -\big(\lambda_0^2 \alpha^2 +\beta^2\lambda_k^2 H_k^2(\xi)\big)\Big]\, d\xi \\
&= (c_5-\lambda_k^2)\int_{S^{n-1}}\Big[\alpha^2 + \beta^2H_k^2(\xi)\Big]\,d\xi +(\lambda_k^2-\lambda_0^2)|S^{n-1}| \alpha^2  \\
&=(c_5-\lambda_k^2)\delta^2+O(\delta^4).
\end{align*}
Thus
\begin{equation*}
\eqref{ball3d}=-c_1^2c_4(c_5-\lambda_k^2)\,\delta^2 + O(\eps\delta^2).
\end{equation*}
Choosing  $k$ large
enough, so that $\lambda_k^2<c_5$ and $\delta>0$ (and therefore $\epsilon$)  small enough, we can guarantee that the sign of \eqref{ball3d} is negative, which means
$$\int_{S^{n-1}}\mathrm{vol}(K\cap\xi^\perp)^n\,d\xi<
\int_{S^{n-1}}\mathrm{vol}(B\cap\xi^\perp)^n\,d\xi.$$
On the other hand, if $k=2$ then $\lambda_2^2 =|S^{n-2}|^2/(n-1)^2>c_5$, which makes the sign of \eqref{ball3d} positive and yields the reversed inequality: $$\int_{S^{n-1}}\mathrm{vol}(K\cap\xi^\perp)^n\,d\xi> \int_{S^{n-1}}\mathrm{vol}(B\cap\xi^\perp)^n\,d\xi.$$

\qed

Next we will solve  the minimizing problem in the class of star
bodies, not necessarily origin-symmetric. For this we need the following lemma. We will use it to show that the constant in Theorem \ref{th:min3d} is optimal.

\begin{Lem}\label{lem:min3d}
Let $C$ be a spherical cap in $S^{n-1}$, $n\ge3$, of the form
\begin{equation*}
C=\{x\in S^{n-1}:\langle x,u\rangle\ge\alpha\},\quad u\in S^{n-1},\quad0<\alpha<1.
\end{equation*}
Then, for each $\lambda\in(0,1)$ and small $\eps>0$, there exists a subset $A$ of $C$ such that $|A|=\lambda|C|$ and for each $\xi\in S^{n-1}$
$$|A\cap\xi^\perp|\le\lambda|C\cap\xi^\perp|+\eps.$$
\end{Lem}

\proof
Let $\delta>0$ be small and $0<\gamma<1$, both of which will be fixed later according to $\eps>0$. Write $C$ as the following union of disjoint narrow strips.
$$C=(A_1\cup A_2\cup\cdots\cup A_N)\cup(B_1\cup B_2\cup\cdots\cup B_N)$$
where $N$ is the smallest integer greater than $(1-\alpha)/\delta$ and, for each $k=1,\ldots,N$,
\begin{align*}
A_k&=\big\{x\in C:\alpha+(k-\gamma)\delta\le\langle x,u\rangle\le \alpha+k\delta\big\}\\
B_k&=\big\{x\in C:\alpha+(k-1)\delta<\langle x,u\rangle<\alpha+(k-\gamma)\delta\big\}.
\end{align*}
Note that  $A_N$ may be empty or have a smaller height than the rest of
$A_k$. $B_N$ cannot be empty, but it may possibly have  a smaller
height than the rest of
$B_k$.

Then we have
$$|A_k|=\int_{\alpha+(k-\gamma)\delta}^{\alpha+k\delta}f_n(t)\,dt \quad\text{and}\quad |B_k|=\int_{\alpha+(k-1)\delta}^{\alpha+(k-\gamma)\delta}f_n(t)\,dt,$$
where $f_n(t)=|S^{n-2}|(1-t^2)^{(n-3)/2}$ if $0\le t\le 1$ and
$f_n(t)=0$ if $t>1$. (This  extension of the function $f_n$ to the
values $t>1$ is dictated by  the possibility that $\alpha + N\delta >1$).

 Since $f_n$ is non-increasing, we have
\begin{align*}
\gamma\delta\cdot f_n\big(\alpha+(k-\gamma)\delta\big) &\ge |A_k| \ge \gamma\delta\cdot f_n\big(\alpha+k\delta\big),\\
(1-\gamma)\delta\cdot f_n\big(\alpha+(k-1)\delta\big) &\ge |B_k| \ge (1-\gamma)\delta\cdot f_n\big(\alpha+(k-\gamma)\delta\big).
\end{align*}
Assuming that $B_{N+1}=\emptyset$, this implies
\begin{equation}\label{eq:altineq}
\frac{|B_k|}{1-\gamma}\ge\frac{|A_k|}{\gamma}\ge\frac{|B_{k+1}|}{1-\gamma},\quad k=1,2,\ldots,N.
\end{equation}

Let $A=A_1\cup A_2\cup\cdots\cup A_N$. Using the left inequality of \eqref{eq:altineq}, we get
\begin{align}
|A|&=\sum_{k=1}^N|A_k| = \sum_{k=1}^N\Big( (1-\gamma)|A_k|+\gamma|A_k|
     \Big)\notag \\
&\le \sum_{k=1}^N\Big(\gamma|B_k|+\gamma|A_k|\Big) =\gamma |C|. \label{eq:upbound}
\end{align}
Also, using the right inequality of \eqref{eq:altineq}, we get
\begin{align}
|A|&=\sum_{k=1}^N|A_k| = \sum_{k=1}^N\Big((1-\gamma)|A_k|+\gamma|A_k|\Big) \notag \\
&\ge \sum_{k=1}^N\Big(\gamma|B_{k+1}|+\gamma|A_k|\Big) =\gamma \big(|C|-|B_1|\big). \label{eq:downbound}
\end{align}
Note that $|A|$ increases from $0$ to $|C|$ as $\gamma$ changes from $0$ to $1$, so we can fix $\gamma=\gamma(\delta)$ with $|A|=\lambda|C|$. Then \eqref{eq:upbound}, \eqref{eq:downbound} imply that $\lambda\le\gamma$ and
\begin{align*}
\gamma &\le \frac{\lambda|C|}{|C|-|B_1|}\le
         \frac{\lambda|C|}{|C|-\delta f_n(\alpha)}
= \frac{\lambda}{1-\delta f_n(\alpha)/|C|}.
%\le\lambda+c_1\delta
\end{align*}
Assuming that $\delta<f_n(\alpha)/(2|C|)$, we get
\begin{equation}\label{eq:lambdaconv}
\gamma\le \lambda\left(1 + \frac{2f_n(\alpha)}{|C|} \delta\right) =
\lambda + c_1 \delta,
\end{equation}
where $c_1=2\lambda f_n(\alpha)/|C|$.

To check the statement for $A\cap\xi^\perp$, we write $\xi\in S^{n-1}$ as $$\xi=u\sqrt{1-s^2}+sv,\quad 0\le s\le1,\quad v\in S^{n-1}\cap u^\perp.$$ If $0\le s<\alpha$, then $C\cap\xi^\perp$ is empty; so assume $s\ge\alpha$.
Note that for $x\in\xi^\perp$,
 $$\langle x,u\rangle=\Big\langle x, u-\langle u,\xi\rangle\xi\Big\rangle=s\langle x,\tilde u\rangle,$$ where
$$\tilde u=\frac{u-\langle u,\xi\rangle\xi}{|u-\langle u,\xi\rangle\xi|}=\frac{u-\langle u,\xi\rangle\xi}{\sqrt{1-\langle u,\xi\rangle^2}}=\frac{u-\langle u,\xi\rangle\xi}s.$$
Let $\tilde\delta=\delta/s$, $\tilde\alpha=\alpha/s$, and $\tilde N$
be the largest integer $k$ such that $A_k\cap\xi^\perp$ is
non-empty. (Note that the definition of $\tilde N$ is slightly
different from that of $N$). If $A_k\cap\xi^\perp$ is empty for all
$k$, then we are done. So we will assume $\tilde N\ge 1$. We have
\begin{align*}
A_k\cap\xi^\perp&=\big\{x\in C\cap\xi^\perp:\alpha+(k-\gamma)\delta\le\langle x,u\rangle\le \alpha+k\delta\big\}\\
&=\big\{x\in C\cap\xi^\perp:\tilde \alpha+(k-\gamma)\tilde\delta\le\langle x,\tilde u\rangle\le \tilde \alpha+k\tilde \delta\big\}
\end{align*}
and, similarly,
\begin{align*}
B_k\cap\xi^\perp&=\big\{x\in C\cap\xi^\perp:\tilde \alpha+(k-1)\tilde
                  \delta<\langle x,\tilde u\rangle<\tilde \alpha+(k-\gamma)\tilde\delta\big\}.
\end{align*}
In the case of $n\ge4$, using the same argument as in
\eqref{eq:altineq} for $A\cap\xi^\perp$, $C\cap\xi^\perp$, $\tilde
\alpha$, $\tilde\delta$, $f_{n-1}$ instead of $A$, $C$, $\alpha$,
$\delta$, $f_n$, we have
$$\frac{|B_k\cap\xi^\perp|}{1-\gamma}\ge\frac{|A_k\cap\xi^\perp|}{\gamma},\quad
k=1,\dots, \tilde N,$$
since the function $f_{n-1}$ is still non-increasing when $n\ge4$. So
$$|A\cap\xi^\perp| \le \sum_{k=1}^{\tilde N}\Big(\gamma|B_k\cap\xi^\perp|+\gamma|A_k\cap\xi^\perp|\Big)\le\gamma|C\cap\xi^\perp|.$$
Using \eqref{eq:lambdaconv} to bound $\gamma$ and letting $c_2=c_1|S^{n-2}|$ yield $$|A\cap\xi^\perp| \le \lambda|C\cap\xi^\perp| + c_2\delta.$$
On the other hand, if $n=3$,  the function $f_{n-1}(t)=2/\sqrt{1-t^2}$ is increasing on $(0,1)$, which implies $$\frac{|A_k\cap\xi^\perp|}{\gamma}\le\frac{|B_{k+1}\cap\xi^\perp|}{1-\gamma},\quad k=1,\ldots,\tilde N-1,$$
and hence $$\sum_{k=1}^{\tilde N-1}|A_k\cap\xi^\perp| \le \sum_{k=1}^{\tilde N-1}\Big(\gamma|B_{k+1}\cap\xi^\perp|+\gamma|A_k\cap\xi^\perp|\Big)\le\gamma|C\cap\xi^\perp|.$$
Thus
\begin{align*}
|A\cap\xi^\perp|&\le \gamma|C\cap\xi^\perp|+|A_{\tilde N}\cap\xi^\perp| \le \gamma|C\cap\xi^\perp|+\int_{1-\tilde\delta}^1\frac{2\,dt}{\sqrt{1-t^2}} \\
&\le \lambda|C\cap\xi^\perp|+c_2\delta+2\arccos(1-\delta/\alpha),
\end{align*}
where we assumed $\delta < \alpha$.

Finally, choosing an appropriate $\delta>0$, satisfying
$c_2\delta+2\arccos(1-\delta/\alpha)\le\eps$ if  $n=3$ and
$c_2\delta\le\eps$ otherwise, yields the desired inequality for
$|A\cap\xi^\perp|$.

\qed

\begin{Thm}\label{th:min3d}
Let $K$ be star-shaped in $\mathbb S^n_+$,  $n\ge3$. Then
$$\int_{S^{n-1}}  \mathrm{vol}(K\cap\xi^\perp)^n\, d\xi\ge c_n\mathrm{vol}(K)^n,$$
where $c_n=2\Gamma(\frac{n+1}2)^n\Gamma(\frac{n}2)^{-n-1}$ is optimal for any fixed volume.
\end{Thm}

\proof
We will start with the following inequality
\begin{equation}\label{phi_ineq}
\frac{\phi_{n-1}(\frac\pi2)}{\phi_n(\frac\pi2)}\phi_n(x)\le
\phi_{n-1}(x), \quad0\le  x\le \frac\pi2,
\end{equation}
where, as before, $\phi_n(x)=\int_0^x\sin^{n-1}t\,dt$.

To prove it, consider the function $$f(x) = \frac{\phi_{n-1}(\frac\pi2)}{\phi_n(\frac\pi2)}\phi_n(x)-
\phi_{n-1}(x).$$
One can see that $f(0) = f(\pi/2)=0$. Moreover, $f'$ changes its sign
exactly once in the interval $(0,\pi/2)$ and $f'(x)>0$ when $x$ is
near $\pi/2$. Thus $f$ is first decreasing and then increasing, which
means that $f(x)<0$ when $x\in (0,\pi/2)$.

Next, using  H\"{o}lder's inequality,  formula (\ref{self-adj}),
and   inequality \eqref{phi_ineq}, we get
\begin{align*}%\label{eq:minHolder}
&\int_{S^{n-1}}\mathrm{vol}(K\cap\xi^\perp)^nd\xi\\
 &\qquad =
\int_{S^{n-1}}  \left(\int_{S^{n-1}\cap\xi^\perp}
                                                   \int_0^{\rho_K(\theta)}
                                                   \sin^{n-2} r \, dr\,
                                                   d\theta \right)^nd\xi\\
&\qquad \ge |S^{n-1}|^{-n+1}  \left(\int_{S^{n-1}}  \int_{S^{n-1}\cap\xi^\perp}
                                                   \phi_{n-1}(\rho_K(\theta))\,
                                                   d\theta \, d\xi \right)^n\\
&\qquad  = |S^{n-1}|^{-n+1} |S^{n-2}|^n  \left(\int_{S^{n-1}}    \phi_{n-1}(\rho_K(\xi))\,
                                                   \, d\xi \right)^n\\
&\qquad \ge |S^{n-1}|^{-n+1} |S^{n-2}|^n
  \left[\frac{\phi_{n-1}(\pi/2)}{\phi_n(\pi/2)}\right]^n  \left(\int_{S^{n-1}}    \phi_{n}(\rho_K(\xi))\,
                                                   \, d\xi \right)^n\\
&\qquad = c_n\mathrm{vol}(K)^n ,
\end{align*}
where $$c_n=
\left[\frac{\phi_{n-1}(\pi/2)}{\phi_n(\pi/2)}\right]^n\frac{|S^{n-2}|^n}{|S^{n-1}|^{n-1}}=\frac{|S^{n-1}|^{n+1}}{|S^n|^n}=\frac{2\Gamma(\frac{n+1}2)^n}{\Gamma(\frac{n}2)^{n+1}}.$$

%To show that the constant is optimal, consider balls $B_r$ of radius
%$r$ that approaches $\pi/2$. Then it is easy to see that $$\lim_{r\to\pi/2}
%\frac{\int_{S^{n-1}}\mathrm{vol}(B_r\cap\xi^\perp)^nd\xi}{\mathrm{vol}(B_r)^n} = c_n.$$

To show that the constant $c_n$ is optimal for {\it any fixed volume} we construct a family of origin-symmetric star-shaped sets of a given volume in $\mathbb S^n_+$, for which the infimum of the quantity $\mathrm{vol}(K)^{-n}\int_{S^{n-1}} \mathrm{vol}(K\cap\xi^\perp)^n$ is equal to $c_n$. Let $0<t<1$ be fixed  and take a small $\eps>0$. Consider a spherical cap $C_\alpha=\{x\in S^{n-1}:\langle x,e_1\rangle\ge\alpha\}$ for small $\alpha>0$ with $|C_\alpha|>\frac{t}2|S^{n-1}|$. Applying Lemma~\ref{lem:min3d} to $C_\alpha$ and $\lambda=\frac{t|S^{n-1}|}{2|C_\alpha|}$, we get $A_t\subset C_\alpha$ such that $|A_t|=\lambda|C_\alpha|$ and $|A_t\cap\xi^\perp|\le\lambda|C_\alpha\cap\xi^\perp|+\eps$ for each $\xi\in S^{n-1}$. Let $K_t$ be the (origin-symmetric) cone with base $A_t\cup(-A_t)$. Then $$\mathrm{vol}(K_t)=2\phi_n(\pi/2)|A_t|=2\lambda\phi_n(\pi/2)|C_\alpha|=t\,\mathrm{vol}(\mathbb S^n_+)$$ and
\begin{align*}
\frac1{\mathrm{vol}(K_t)^n}\int_{S^{n-1}}\mathrm{vol}(K_t\cap\xi^\perp)^nd\xi = \int_{S^{n-1}}\Big(\frac{2\phi_{n-1}(\frac\pi2)|A_t\cap\xi^\perp|}{2\lambda\phi_n(\frac\pi2)|C_\alpha|}\Big)^nd\xi \\
\le\frac{\phi_{n-1}^n(\frac\pi2)}{\phi_n^n(\frac\pi2)} \int_{S^{n-1}}\Big(\frac{\lambda|C_\alpha\cap\xi^\perp|+\eps}{\lambda|C_\alpha|}\Big)^nd\xi.
\end{align*}
As $\alpha\to0$ and $\eps\to0$, the right hand side of the above inequality approaches the constant $c_n$.

\qed

\begin{Rmk} Even though the constant $c_n$ in Theorem~\ref{th:min3d} is optimal for star-shaped sets of any fixed volume,  the inequality in Theorem~\ref{th:min3d} is strict  for most of them.   Equality occurs only in special cases, namely,
if and only if $K$ is a cone whose base
$B\subset S^{n-1}$ is one of the following types:
\begin{enumerate}
\item $B$ or its complement in $S^{n-1}$ is of measure zero.
\item The intersection of  $B$ and the reflection of its complement
  with respect to the origin is of measure zero.
\end{enumerate}
\end{Rmk}

\proof
Equality in \eqref{phi_ineq} holds only when $x=0$ or $x=\pi/2$. Thus,
$\rho_K$ can only take two values: $0$ and $\pi/2$. Furthermore,
equality in  H\"{o}lder's inequality occurs when
 $ \int_{S^{n-1}\cap\xi^\perp}  \phi_{n-1}(\rho_K(\theta))\,
 d\theta$
is a constant for almost all $\xi\in S^{n-1}$.
This means that $\phi_{n-1}(\rho_K(\theta))+
\phi_{n-1}(\rho_K(-\theta))$ is also a constant for almost all $\theta \in S^{n-1}$. The
value of this constant can only be $0$,  $\phi_{n-1}(\pi/2)$,
$2\phi_{n-1}(\pi/2)$. Clearly, in the first case we get
$\rho_K(\theta) = 0$, in the second case   $\rho_K(\theta) =
0$ if and only if $\rho_K(-\theta) = \pi/2$, in the third case
$\rho_K(\theta) = \pi/2$ (up to sets of measure zero).

\qed

Using Lemma \ref{lem:min3d} we can also show that the infimum of  $\int_{S^{n-1}}  |K\cap\xi^\perp|^n\, d\xi$ in the class of (origin-symmetric) star bodies in $\mathbb R^n$ or $\mathbb H^n$ of any fixed volume is zero.

\begin{Thm}
Let $\mathbb M^n=\mathbb R^n$ or $\mathbb H^n$ for $n\ge3$ and $\alpha>0$.
For each $\eta>0$ there is an origin-symmetric star body $K\subset\mathbb M^n$ of volume $\alpha$ such that
$$\int_{S^{n-1}}  \mathrm{vol}(K\cap\xi^\perp)^n\, d\xi\le \eta.$$
\end{Thm}

\proof
Let $\phi_n(x)=\int_0^x s(t)^{n-1}\,dt,$ where $s(t)=t$ if $\mathbb M^n=\mathbb R^n$, and $s(t)=\sinh t$ if $\mathbb M^n=\mathbb H^n$. In either case  we have
$$\lim_{x\to\infty}\phi_n(x)=\infty, \quad \lim_{x\to\infty}\frac{\phi_{n-1}(x)}{\phi_n(x)}=\lim_{x\to\infty}\frac{\phi_{n-1}'(x)}{\phi_n'(x)}=0.$$
Fix a spherical cap $C\subset S^{n-1}$ of measure less than $\frac12|S^{n-1}|$. Let $r$ be large enough so that $\phi_n(r)>\frac{\alpha}{2|C|}$. Applying Lemma \ref{lem:min3d} for the spherical cap $C$, $\lambda=\frac{\alpha}{2\phi_n(r)|C|}$, and $\eps=\frac1{\phi_n(r)}$, we can construct a set $A\subset S^{n-1}$ such that $|A|=\lambda|C|$ and for each $\xi\in S^{n-1}$
\begin{align*}
|A\cap\xi^\perp|&\le\lambda|C\cap\xi^\perp|+\eps=\frac1{\phi_n(r)}\cdot\frac{\alpha|C\cap\xi^\perp|+2|C|}{2|C|}.
\end{align*}

Define the origin-symmetric star-shaped set $K\subset\mathbb M^n$ by
\begin{equation*}
\rho_K(x)=
\begin{cases}
r,& \text{if } x\in A\cup(-A),\\
0,& \text{otherwise.}
\end{cases}
\end{equation*}
Then $$\mathrm{vol}(K)=\int_{S^{n-1}}\phi_n(\rho_K(x))\,dx=2\phi_n(r)|A|=2\phi_n(r)\lambda|C|=\alpha$$
and
\begin{align*}
\int_{S^{n-1}}\mathrm{vol}(K\cap\xi^\perp)^n\,d\xi &=\int_{S^{n-1}}\left(2\phi_{n-1}(r)|A\cap\xi^\perp|\right)^n\, d\xi \\
&\le\left(\frac{\phi_{n-1}(r)}{\phi_n(r)}\right)^n\int_{S^{n-1}}\left(\frac{\alpha|C\cap\xi^\perp|+2|C|}{|C|}\right)^n\, d\xi.
\end{align*}
Since the integral in the right-hand side is a constant and $\frac{\phi_{n-1}(r)}{\phi_n(r)}\to0$ as $r\to\infty$, choosing $r$ large enough, so that the right hand side is less than $\eta$, yields the desired inequality. The statement is now proved for star-shaped sets. By an approximation the statement extends to  star bodies.

\qed

\noindent{\bf Question.} Is it true that the maximizers in the class of origin-symmetric convex
bodies in $\mathbb S^n_+$, $n\ge 3$, of a fixed volume are the lunes?

\section{Busemann's intersection inequality in a general measure
  space}

In what follows,  $\mu$ will be a measure on $\mathbb M^n=\mathbb R^n$
or $\mathbb H^n$  with density $f$. If $K$ is a star body in $\mathbb
M^n$, then $$\mu(K)  = \int_K f(x)\, d\mathrm{vol}(x)$$ and $$\mu(K\cap \xi^\perp)
= \int_{K\cap \xi^\perp}
f(x)\, d\mathrm{vol}(x),$$  where  $d\mathrm{vol}$ is the volume element on $\mathbb M^n$
or the submanifold $\xi^\perp$.

\begin{Thm}\label{General_m}
Let $K$ be a star body in $\mathbb M^n=\mathbb R^n$ or $\mathbb H^n$ 
and $\mu$ be a measure on $\mathbb M^n$, whose density $f$ is radially symmetric and decreasing. Then
$$\int_{S^{n-1}}\mu(K\cap\xi^\perp)^n\,d\sigma_{n-1}(\xi) \le \left[\Psi(\mu(K))\right]^{n-1}$$
with equality if and only if $K$ is a ball centered at the
origin. Here, $d\sigma_{n-1}$ is the spherical volume element
normalized so that $\sigma_{n-1}(S^{n-1}) =1$, and $\Psi=\psi_{n-1}^{\frac{n}{n-1}}\circ\psi_n^{-1}$ is a
positive increasing concave function on $[0,\infty)$,  where $\psi_n(x)$ denotes the
measure (with respect to $\mu$) of the centered ball
in $\mathbb M^n$ of radius $x\ge0$.
\end{Thm}

\proof
The function $\psi_n$ can be written
as $$\psi_n(x)= |S^{n-1}| \int_0^xf(t)s(t)^{n-1}dt,$$ where $f$ is decreasing on
$[0,\infty)$, $s(t)=   t$ if $\mathbb M^n=\mathbb
R^n$ and $s(t)= \sinh t$ if $\mathbb M^n=\mathbb H^n$. Note that in
either case  $s(t)$ is a convex function on $[0,\infty)$. First we show that $\Psi$ is concave. Differentiating $\Psi(\psi_n(x))=[\psi_{n-1}(x)]^{\frac{n}{n-1}}$ we have
$$\Psi'(\psi_n(x))=\frac{|S^{n-2}|}{|S^{n-1}|}\frac{n}{n-1}[\psi_{n-1}(x)]^{\frac1{n-1}}/s(x)$$
and
$$\Psi''(\psi_n(x))=\frac{n
  |S^{n-2}|[\psi_{n-1}(x)]^{(2-n)/(n-1)}}{(n-1) |S^{n-1}|^2
  f(x)s(x)^{n+1}}\left[\frac{|S^{n-2}|f(x)s(x)^{n-1}}{n-1}-\psi_{n-1}(x)s'(x)\right].$$ Since $f$ is decreasing and $s$ is convex, we have
\begin{align*}
\psi_{n-1}(x)s'(x)&=s'(x)|S^{n-2}| \int_0^xf(t)s(t)^{n-2}dt \\
&> f(x)|S^{n-2}| \int_0^xs(t)^{n-2}s'(t)dt=\frac{|S^{n-2}|f(x)s(x)^{n-1}}{n-1}.
\end{align*}
This implies $\Psi''<0$, so $\Psi=\psi_{n-1}^{\frac{n}{n-1}}\circ\psi_n^{-1}$ is concave on $(0,\infty)$.

Let $L$ be the star body  in $\mathbb R^n$ with the radial function
$$\rho_L(u)=\left[\frac{\psi_{n-1}(\rho_K(u))}{|B_2^{n-1}|}\right]^{\frac1{n-1}},\quad u\in S^{n-1}.$$
Then, for each $\xi\in S^{n-1}$,
\begin{align*}
\mu(K\cap\xi^\perp)&=\int_{S^{n-1}\cap\xi^\perp}\int_0^{\rho_K(u)}f(t)s(t)^{n-2}\,dt\,du\\
&=\frac{1}{|S^{n-2}|}\int_{S^{n-1}\cap\xi^\perp}\psi_{n-1}(\rho_K(u))\,du\\
&=\frac1{n-1}\int_{S^{n-1}\cap\xi^\perp}\rho_L^{n-1}(u)\,du =|L\cap\xi^\perp|^n.
\end{align*}
Busemann's intersection inequality for $L\subset\mathbb R^n$ implies
\begin{equation}\label{eq:busemannineq}
\int_{S^{n-1}}|L\cap\xi^\perp|^n\,d\sigma_{n-1}(\xi)\le c_1|L|^{n-1},
\end{equation}
where $c_1=|B_2^{n-1}|^n/|B_2^n|^{n-1}$.
The concavity of $\Psi$ implies
\begin{align}\label{eq:jensenineq}
|L|&=\frac1n\int_{S^{n-1}}\rho_L^n(u)\,du=\frac1n\int_{S^{n-1}}\left[\frac{\psi_{n-1}(\rho_K(u))}{|B_2^{n-1}|}\right]^{\frac{n}{n-1}}\,du \notag \\
&=c_2\int_{S^{n-1}}(\psi_{n-1}^{\frac{n}{n-1}}\circ\psi_n^{-1})(\psi_n(\rho_K(u)))\, d\sigma_{n-1}(u)\notag \\
&\le c_2\Psi\left(\int_{S^{n-1}}\psi_n(\rho_K(u))\, d\sigma_{n-1}(u)\right)=c_2\Psi(\mu(K)),
\end{align}
with $c_2=|B_2^n|/|B_2^{n-1}|^{n/(n-1)}$. Therefore
$$\int_{S^{n-1}}\mu(K\cap\xi^\perp)^n\,d\sigma_{n-1}(\xi)\le c_1\left[c_2\Psi(\mu(K))\right]^{n-1}=\left[\Psi(\mu(K))\right]^{n-1}.$$
Moreover,  the equality case in   inequalities
\eqref{eq:busemannineq} and \eqref{eq:jensenineq} holds if and only if
$K$ is a centered ball (we use that $\Psi$ is strictly concave).

\qed

\begin{Rmk} An important  particular case of the previous theorem is
  the case when the measure $\mu$ is the standard Gaussian measure, i.e.,
$$\mu(K) = (2\pi)^{-n/2} \int_K e^{-|x|^2/2} dx,$$  $$  \mu(K\cap \xi^\perp) = (2\pi)^{(-n+1)/2} \int_{K\cap\xi^\perp}  e^{-|x|^2/2} dx.$$
Note that one has to adjust the coefficient in Theorem  \ref{General_m}  since the densities on $\mathbb R^n$ and $\mathbb R^{n-1}$ differ by a factor of $\sqrt{2\pi}$.

\end{Rmk}

\end{document}